\documentclass[11pt]{article}
\usepackage[square,semicolon,authoryear]{natbib}

\usepackage{amsthm}
\usepackage{amsfonts}
\usepackage{amssymb}
\usepackage{graphicx}
\usepackage{caption}
\usepackage{amsmath}
\usepackage{lscape}
\usepackage{epstopdf}
\usepackage{subfigure}
\usepackage{fullpage}
\usepackage{tikz}
\newtheorem{Theorem}{Theorem}[section]

\newtheorem{Corollary}[Theorem]{Corollary}
\newtheorem{Lemma}[Theorem]{Lemma}
\newtheorem{Example}[Theorem]{Example}

\newtheorem{Remark}[Theorem]{Remark}

\begin{document}

\title{\bf Optimal Redundancy Allocation in  Coherent Systems with Heterogeneous Dependent Components\\
$To~ appear~ in ~J. ~Applied ~Probability$ }

\author{\textbf{Maryam Kelkinnama\footnote{Department of Mathematical Sciences, Isfahan University of Technology, Isfahan 84156-83111, Iran,
(Email: m.kelkinnama@iut.ac.ir )} \ \ and\ \ Majid Asadi\footnote{{Department of Statistics, University of Isfahan, Isfahan 81744, Iran, \& School of  Mathematics, Institute of Research in Fundamental Sciences (IPM), P.O Box 19395-5746, Tehran, Iran,} (E-mail: m.asadi@sci.ui.ac.ir)}, }\\}
\date{}
\maketitle

\begin{abstract}
This paper is concerned with the optimal number of redundant allocation to $n$-component coherent systems consist of heterogeneous dependent components. We assume that the system is built of $L$ groups of different components, $L\geq 1$, where there are $n_i$ components in group $i$, and $\sum_{i=1}^{L}n_i=n$. The problem of interest is to allocate $v_i$ active redundant components to each component of type $i$, $i=1,\dots, L$. To get the optimal values of $v_i$, we propose two cost-based criteria. One of them is introduced based on the costs of renewing the failed components and the costs of refreshing the alive ones at the system failure time. The other criterion is proposed based on the costs of replacing the system at its failure time or at a predetermined time $\tau$, whichever occurs first. The expressions for the proposed functions are derived using the mixture representation of the system reliability function based on the notion of survival signature. We assume that a given copula function models the dependency structure between the components. In the particular case that the system is a series-parallel structure, we provide the formulas for the proposed cost-based functions. The results are discussed numerically for some specific coherent systems.

\end{abstract}

\date{}
\maketitle

\section{Introduction}
\subsection{Motivation and related literature}
\quad In reliability engineering and system security, one of the most useful methods for enhancing the reliability characteristics of a system is to allocate redundant components to the system.
The redundancy can be performed at the component level or the system level. In the former case, some redundant components are connected to each component, while in the latter one, the original coherent system fastens to some copies of itself. In a commonly used type of redundancy, called active redundancy, the original component and the redundant ones work simultaneously as parallel. In this case, the lifetime of the resulted parallel subsystem equals the maximum lifetime between the connected components. This strategy is mostly applied when the replacement of the components during the operation time of the system is impossible.
As redundancy allocation is a widely used method for improving the performance of the products, numerous researchers have paid attention to develop the theories and applications of this subject.
For example,  \cite{li-ding2010}  investigated the allocation of active redundancies to a $k$-out-of-$n$ system in which the lifetimes of independent components are stochastically ordered. 
\cite{you2016} studied $k$-out-of-$n$ redundant systems with dependent components.
{\cite{Eryilmaz-Ucum} determined the optimal number of spare components for a  weighted $k$-out-of-$n$.} 
 \cite{kavlak} investigated the reliability and the mean residual life  functions of coherent systems with active redundancies at the component and system levels.
 \cite{zhang} investigated the optimal allocation of active redundancies for weighted $k$-out-of-$n$ systems. \cite{amini} compared the component redundancy versus system redundancy for coherent systems with dependent and identically distributed components.  \cite{fangli2016} studied the allocation of one active redundancy to coherent systems consisting of heterogeneous and statistically dependent components. Utilizing the minimal path decomposition, they proposed a necessary and sufficient condition identifying a better allocating strategy from two candidates.  \cite{fangli2017} investigated allocating multiple matched active redundant components to coherent systems.  \cite{fangli2018} studied the coherent systems with one active redundancy, using the minimal cut decomposition of the system. { \cite{Torradoetal2021}  studied the redundancy allocation for a coherent system formed by modules, under  different settings related to dependency  and distribution of components. They stochastically compared the redundancies at the components level versus redundancies at the modules level.}
{ \cite{Torrado} considered a coherent system having possibly dependent subsystems in which the components are connected in parallel or in series. 
It is assumed that a number of possibly dependent components in each subsystem  are randomly selected from a heterogeneous population. 
The cited author stochastically compared such systems with different numbers of components, based on majorization orders and, determined the optimal numbers of components in each subsystem such that the system  reliability is maximized. Particularly, she examined the results for series-parallel systems. 
} The redundancy allocation  in a series-parallel system has also been considered by some authors, among which we refer to \cite{soltani-O-2014}, \cite{karimi}, and \cite{fang2020}. 


\quad It is worth noting that the redundant components can be added to the system as inactive (cold and warm standby) components.  The systems with cold and warm standby redundancy have also investigated in reliability literature; see, for example,  \cite{eryi-cold},   \cite{fink2018}, \cite{shen}, and \cite{asadi}. 

 \subsection{Survival signatures of coherent systems}\label{surv}
\quad The first main step to assess the reliability and stochastic characteristics of an $n$-component system is to get knowledge about the structure-function of the system as well as the probability distribution of component lifetimes.
In this regard, a useful concept for assessing the reliability of the system through the reliability of its components is the notion of {\it survival signature}.  This concept is particularly significant for describing the structures of coherent systems with multiple types of components.
  Consider an $n$-component coherent system consisting of $L$ different types, such that there are $n_i$ components from $i$th type, $i=1,\ldots, L$,  and $\sum_{i=1}^L n_i=n$.
The reliability function of the system, at any time $t$, can be represented as follows
\begin{align}\label{Fbaroriginal}
\bar{F}_T(t)=\sum_{l_1=0}^{n_1}\ldots\sum_{l_L=0}^{n_L}\Phi(l_1,\ldots, l_L)P(C_1(t)=l_1, \ldots, C_L(t)=l_L),
\end{align}
where  $C_i(t)$ denotes the number of components of type $i$ working at time $t$, and $\Phi$ is called the survival signature and represents \lq\lq the probability that the system is working when exactly $l_i$ components of type $i$ is working\rq\rq, see \cite{coolen}.

Suppose that  the  lifetimes of the components of the same type are exchangeable dependent, and the lifetimes of the components of different types are dependent.
Commonly, the dependency structure is modeled using a survival copula. In other words, if  $T_j^{(i)}$ denotes the lifetime of $j$th component from type $i$, $j=1, \ldots, n_i$, $i=1, \ldots, L$, and $\bar{F}_i,~i=1, \ldots, L$ denotes the common reliability function for the components of the $i$th type, then there is a survival copula  $\hat{C}$ such that the joint reliability  of $T_j^{(i)}$'s can be written as
 \begin{align}\label{copula}
&P( T_1^{(1)}>t_1^{(1)}, \ldots,  T_{n_1}^{(1)}>t_{n_1}^{(1)}, \ldots,  T_1^{(L)}>t_L^{(L)}, \ldots,  T_{n_L}^{(L)}>t_{n_L}^{(L)})\nonumber\\
 &\hspace{2cm}= \hat{C}\big(\bar{F}_1(t_1^{(1)}), \ldots, \bar{F}_1(t_{n_1}^{(1)}), \ldots, \bar{F}_L(t_1^{(L)}), \ldots, \bar{F}_L(t_{n_L}^{(L)})\big),
 \end{align}
see, for example, \cite{navarro2}, \cite{navarro1}, and \cite{fangli2018}.  In this case, it can be shown that
 \begin{align}\label{Fbar-dep}
\bar{F}_T(t)&={\sum_{l_1=0}^{n_1}\ldots\sum_{l_L=0}^{n_L}}{\sum_{i_1=0}^{n_1-l_1}\ldots\sum_{i_L=0}^{n_L-l_L}}(-1)^{i_1+\ldots+ i_L}\binom{n_1}{l_1}\ldots \binom{n_L}{l_L}\binom{n_1-l_1}{i_1}\ldots \binom{n_L-l_L}{i_L}\nonumber\\
&\times\Phi(l_1,\ldots, l_L) \hat{C}(\underbrace{\bar{F}_1(t)}_{i_1+l_1},\underbrace{1}_{n_1-(i_1+l_1)},\ldots,\underbrace{\bar{F}_L(t)}_{i_L+l_L},\underbrace{1}_{n_L-(i_L+l_L) }),
\end{align}
where  $\underbrace{u}_{m}$ means the $m$ repetitions of $u$,
see \cite{eryilmaz-importance-2018, eriylmaz-dependent-2018}.
If the components of the system are independent, then the representation \eqref{Fbaroriginal} is converted to the following expression
\begin{align}\label{ind}
\bar{F}_T(t)=\sum_{l_1=0}^{n_1}\ldots\sum_{l_L=0}^{n_L}\Phi(l_1,\ldots, l_L)\prod_{i=1}^L \binom{n_i}{l_i}[\bar{F}_i(t)]^{l_i} [{F}_i(t)]^{n_i-
l_i}.
\end{align}

Many authors have considered the reliability properties of a coherent system with multi-type components based on the survival signature and for various applications. Among them, we mention to recent papers, \cite{feng2016}, \cite{saman}, \cite{eryilmaz-importance-2018, eriylmaz-dependent-2018}.
  \cite{heu} used
the notion of survival signature for the formulation of the reliability-redundancy allocation problem. They considered the objective function to maximize the system reliability under some constraints.  

\quad \cite{zarezadeh2019}  studied the reliability and preventive maintenance of the coherent systems with multi-type components whose components are subject to failure according to multiple external shocks. \cite{hashemi} proposed two maintenance strategies for optimal preservation of coherent systems consisting of independent multi-type components.

\subsection{Contributions of this paper}
\quad
This paper aims to study the optimal number of redundancy allocation to $n$-component coherent systems consist of different components. 
It is assumed that the components of the system are dependent, where a given copula function models the dependency structure.
 We are interested in allocating $v_i$ active redundant components to each component of type $i$, under the 
 constraint on the number of existing spare components. To get the optimal number of $v_i$'s,  we propose two cost-based functions. 
More precisely,  the contributions of the paper are as follows.
\begin{itemize}
   \item We propose a mean cost rate function in terms of the costs of renewing the failed components and the costs of  refreshing the alive components at the time of the system failure. Then, we find the optimal number of redundant components, $v_i$'s, to be added to each component of type $i$, such that the proposed cost function is minimized.

       \item We introduce a mean cost rate function, relevant to an age replacement policy, in terms of the costs of renewing (refreshing) the failed (alive) components at the failure time of the system or at a predetermined time $\tau$, whichever occurs first. Then, the optimal values of $v_i$'s are obtained, such that the suggested cost-based function achieves its minimum value.
           
\item  In the particular considerable case that the system is a series-parallel system, we provide the formulas for the proposed mean cost rate functions.
Then we investigate the optimal number of the components for each parallel subsystem such that the proposed  functions are minimized
 \end{itemize}
   The derivations of the paper are extensions of the results of \cite{Eryilmaz-main-2018}, who investigated  the optimal number of components in the case that  structure function is $k$-out-of-$n$ with independent components.

\subsection{Organizations of the paper }
\quad The remaining of the paper is arranged as follows. In Section 2, under the settings of subsection \ref{surv},  we present the formulation of the system reliability function (\ref{Fbar-dep}) in the case that $v_i$ components are added as active redundant to each component of type $i$, $i=1,\dots, L$.  Then, utilizing this formulation, a mean cost rate function is introduced at the time of the system failure.
 Next, a mean cost rate function is established based on the costs of replacing the system 
 at its failure time or at a predetermined time $\tau$, whichever occurs first. The expressions for the proposed mean cost rate functions are derived in terms of the reliability function (\ref{Fbar-dep}). 
{Some examples of coherent systems are presented to illustrate the applications of the proposed approaches; a 6-component system consisting of two types of dependent components, and an 8-component system composed of three types of components that are independent.} The optimal number of redundant components, based on the proposed cost-based functions, are discussed for each system numerically.

 \quad Section 3 is devoted to the particular case that the system is a series-parallel system. We provide the formulas for the proposed mean cost rate functions in Section 2 for such systems. Then, we investigate the optimal number of the components for each parallel subsystem such that the proposed cost functions are minimized. The results of this section are numerically illustrated for a series-parallel system consisting of three parallel subsystems connected in series. Some concluding remarks in Section 4 finalize the paper. The details of the proofs are given in the Appendix.
\section{Optimal number of redundant components}
\quad We consider an $n$-component coherent system consisting of multiple types of components with following description.
{ The system is built up of $L$ types of components, $L\geq 1$, such that there are $n_i$  components of type $i$  and $\sum_{i=1}^{L}n_i=n$. We assume that the common reliability function of the components of type $i$ is $\bar{F}_i(.)$, $i=1,2,\dots,L$.
The  lifetimes of the components of the same type are exchangeable dependent, and the lifetimes of the components of different types are dependent. The assumed  dependence structure is modeled by a survival copula given in \eqref{copula}. 
 To increase the reliability of the system, we desire to add $v_i$ active redundancies to each component of type $i$, $i=1,\ldots, L$.
 Each original component in the system and its redundant components are assumed to be  independent and identically distributed  (i.i.d.). 
Let $T_R$ denote the lifetime of the system incorporated by redundant components. Because  an original component  and its redundant ones make a parallel subsystem, one can easily show that the reliability function of $T_R$ at time $t$ can be represented as follows.
}
\begin{align*}
\bar{F}_{T_R}(t)
&={\sum_{l_1=0}^{n_1}\ldots\sum_{l_L=0}^{n_L}}{\sum_{i_1=0}^{n_1-l_1}\ldots\sum_{i_L=0}^{n_L-l_L}}(-1)^{i_1+\ldots+ i_L}\binom{n_1}{l_1}\ldots \binom{n_L}{l_L}\binom{n_1-l_1}{i_1}\ldots \binom{n_L-l_L}{i_L} \nonumber\\
&\times\Phi(l_1,\ldots, l_L)\hat{C}(\underbrace{1-{F}_1^{v_1+1}(t)}_{i_1+l_1}, \underbrace{1}_{n_1-(i_1+l_1)},\ldots,\underbrace{1-{F}_L^{v_L+1}(t)}_{i_L+l_L},\underbrace{1}_{n_L-(i_L+l_L) }).
\end{align*}

In the case of independence of all components, this representation reduces to
\begin{align*}
\bar{F}_{T_R}(t)=\sum_{l_1=0}^{n_1}\ldots\sum_{l_L=0}^{n_L}\Phi(l_1,\ldots, l_L)\prod_{i=1}^L \binom{n_i}{l_i}[1-{F}_i^{v_i+1}(t)]^{l_i}
[{F}_i^{v_i+1}(t)]^{n_i-l_i}.
\end{align*}

\quad The problem of interest in this redundancy strategy is to determine the optimal number of spares allocated to each component. In this paper, our approach is to find  $v$'s based on the minimization of a kind of cost criterion. In this regard, we set up two mean cost rate functions to obtain the optimal number of redundant components. One of them  is imposed based on the cost of the system failure, which depends on the number of failed components when a system failure occurs. The other one is defined based on an age replacement policy. In the next subsections, we describe these two functions with details.

	\begin{Remark}
		\em	Although the system considered above is described in the general case that the component lifetimes of the same type are exchangeable dependent, and the lifetimes of the components of different types are dependent, in allocating the redundant components we assumed that  in the constructed parallel subsystem the components are independent and identically distributed. This assumption seems to be a restriction in some practical cases, but it should be noted that if we  drop the assumption of i.i.d.  for the redundant components, the computation of the system reliability,  would be a challenging problem and potentially involves complex calculations. We believe that, considering the problem of optimal redundancy under i.i.d. components  in each subsystem, as considered in this paper, could be a first step toward solving the more general cases (see also \cite{sam}, pp. 76-77).
\end{Remark}
\subsection{Cost function at the system failure}\label{cost}
\quad Suppose that the system starts working at $t=0$ and fails at a random time after $t=0$. Assume that when the system fails, we have a cost $c_i$ for each failed component of type $i$ to replace it by a new one and a cost $c_i^*$ for each unfailed component for refreshing it so that it becomes as good as new, where we assume that $c_i\geq c^*_i$, $i=1, \ldots, L$. Furthermore, we assume that $c^{**}$ denotes the fixed overall cost for system failure. With $T_R$ as the lifetime of the system after redundancy, let the random variable $X_i(T_R)$ denote the number of failed components of type $i$ at the time of system failure, $i=1,\dots,L$. Then, the mean cost rate function for a failed system is defined as
\begin{align}\label{cost1}
Cost_1(\mathbf{v})&=\frac{\sum_{i=1}^Lc_iE(X_i(T_R))+\sum_{i=1}^L c_i^*E(n_i(v_i+1)-X_i(T_R))+c^{**}}{E(T_R)}
\end{align}
where $\mathbf{v}=(v_1,\ldots, v_L)$. The numerator is the expected cost of the system failure, and the denominator is the mean time to failure ($MTTF$) of the system, hence  $Cost_1$ becomes the mean cost per unit of time. Note that in the system after redundancy, { there are altogether $n_i(v_i+1)$ components} of type $i$, $i=1,2, \ldots, L$. The  relation \eqref{cost1} can be rewritten in terms of the lifetime of the original system without any redundancy, $T$, as
\begin{align}\label{poi}
Cost_1(\mathbf{v})&=\frac{\sum_{i=1}^Lc_i(v_i+1)E(X_i(T))+\sum_{i=1}^L c_i^*(v_i+1)E(n_i-X_i(T))+c^{**}}{E(T_R)}\nonumber\\
&=\frac{\sum_{i=1}^L(c_i-c_i^*)(v_i+1)E(X_i(T))+\sum_{i=1}^L c_i^*(v_i+1)n_i+c^{**}}{E(T_R)}.
\end{align}

{\begin{Lemma}\label{ui}
The quantity  $E(X_i(T))$ in (\ref{poi}) can be expressed  as follows.
\begin{small}
\begin{align*}
&E(X_i(T))\\
&=n_i\int_{0}^{\infty} \lim_{\delta\rightarrow 0}\frac{1}{\delta}\sum_{m_1=0}^{n_1}...\sum_{m_i=0}^{n_i-1}...\sum_{m_L=0}^{n_L}\Phi(m_1,...,m_{i-1},m_i+1,m_{i+1},..., m_L) \binom{n_1}{m_1}... \binom{n_i-1}{m_i}...\binom{n_L}{m_L}A_{\mathbf{m}}^{(i)}(t,\delta)dt
\end{align*}
\end{small}
where
\begin{small}
\begin{align}\label{Adelta}
A_{\mathbf{m}}^{(i)}(t,\delta)
&=P\Big({T_1^{(1)}>t},\ldots,{T_{m_1}^{(1)}>t},{T_{m_1+1}^{(1)}\leq t},\ldots, {T_{n_1}^{(1)}\leq t}, \nonumber\\
&\hspace{1cm}\ldots, t<T_1^{(i)}\leq t+\delta, {T_2^{(i)}>t},\ldots,{T_{m_i+1}^{(i)}>t},{T_{m_i+2}^{(i)}\leq t},\ldots, {T_{n_i}^{(i)}\leq t}, \nonumber\\
&\hspace{1cm}\ldots, {T_1^{(L)}>t},\ldots,{T_{m_L}^{(L)}>t},{T_{m_L+1}^{(L)}\leq t},\ldots, {T_{n_L}^{(L)}\leq t}\Big).
\end{align}
\end{small}
\end{Lemma}
\begin{proof}
\begin{align*}
E(X_i(T))&=E(\sum_{j=1}^{n_i}I(T_j^{(i)}\leq T))=\sum_{j=1}^{n_i}P(T_j^{(i)}\leq T)=n_iP(T_1^{(i)}\leq T)\\
&=n_i\int_{0}^{\infty} \lim_{\delta\rightarrow 0}\frac{P(T>t,t<T_1^{(i)}\leq t+\delta)}{\delta}dt,
\end{align*}
where the third equality follows from the exchangeability of the components of type $i$, $i=1, \dots, L$.
By conditioning on the number of live components of each type,  we obtain
\begin{small}
\begin{align}\label{Atdelat}
&P(T>t,t<T_1^{(i)}\leq t+\delta)\nonumber\\
&=\sum_{m_1=0}^{n_1}\ldots\sum_{m_i=0}^{n_i-1}\ldots\sum_{m_L=0}^{n_L}P(T>t,t<T_1^{(i)}\leq t+\delta, C_1(t)=m_1,\ldots, C_i(t)=m_i,\ldots, C_L(t)=m_L)\nonumber\\
&=\sum_{m_1=0}^{n_1}\ldots\sum_{m_i=0}^{n_i-1}\ldots\sum_{m_L=0}^{n_L}\Phi(m_1,\ldots,m_{i-1},m_i+1,m_{i+1},\ldots, m_L) \binom{n_1}{m_1}\ldots \binom{n_i-1}{m_i}\ldots\binom{n_L}{m_L}A_{\mathbf{m}}^{(i)}(t,\delta)
\end{align}
\end{small}
The last equality in \eqref{Atdelat} holds because  the  components of the same type have a common failure time distribution.
\end{proof}}
In the following theorem, \eqref{Adelta} is represented based on the survival copula of the components lifetimes.
\begin{Theorem}\label{Am}
Using the inclusion-exclusion rule, $A_{\mathbf{m}}^{(i)}(t,\delta)$ can be represented as the following expression
\begin{small}
\begin{align*}
A_{\mathbf{m}}^{(i)}(t, \delta)&={\sum_{j_1=0}^{n_1-m_1}\ldots\sum_{j_i=0}^{n_i-m_i-1}\ldots\sum_{j_L=0}^{n_L-m_L}}(-1)^{j_1+\ldots + j_L}\binom{n_1-m_1}{j_1}\ldots \binom{n_i-m_i-1}{j_i}\ldots\binom{n_L-m_L}{j_L}\\
&\times \Big[\hat{C}(\underbrace{\bar{F}_1(t)}_{m_1+j_1},\underbrace{1}_{n_1-(m_1+j_1)},\ldots,\underbrace{\bar{F}_i(t)}_{m_i+j_i+1},\underbrace{1}_{n_i-(m_i+j_i+1)},\ldots,\underbrace{\bar{F}_L(t)}_{m_L+j_L}, \underbrace{1}_{n_L-(m_L+j_L)})\\
&-\hat{C}(\underbrace{\bar{F}_1(t)}_{m_1+j_1},\underbrace{1}_{n_1-(m_1+j_1)},\ldots,\underbrace{\bar{F}_i(t)}_{m_i+j_i},\bar{F}_i(t+\delta),\underbrace{1}_{n_i-(m_i+j_i+1)}, \ldots,\underbrace{\bar{F}_L(t)}_{m_L+j_L},\underbrace{1}_{n_L-(m_L+j_L)})\Big].
\end{align*}
\end{small}
\end{Theorem}
\begin{proof}
See the Appendix.
\end{proof}

Note that, in the  particular case of independence of all components, we get
\begin{small}
\begin{align}\label{1}
&E(X_i(T))\nonumber\\
&=n_i\sum_{m_1=0}^{n_1}\ldots\sum_{m_i=0}^{n_i-1}\ldots\sum_{m_L=0}^{n_L}\Phi(m_1,\ldots,m_{i-1},m_i+1,m_{i+1},\ldots, m_L)\binom{n_1}{m_1}\ldots \binom{n_i-1}{m_i}\ldots\binom{n_L}{m_L}\nonumber\\
&\times\int_{0}^{\infty}\bar{F}_1^{m_1}(t)F_1^{n_1-m_1}(t)\ldots\bar{F}_i^{m_i}(t)F_i^{n_i-m_i-1}(t)\ldots \bar{F}_L^{m_L}(t)F_L^{n_L-m_L}(t)dF_i(t).
\end{align}
\end{small}
\quad In order to minimize the  mean cost rate function $Cost_1(\mathbf{v})$, we impose the constraint that there are at most $M_i$ components of type $i$  as spares, $i=1,\dots,L$. This means that the number of the components that can be connected in parallel at the $i$th group satisfies in the inequality   $n_iv_i\leq M_i$, $i=1,\ldots, L$. To determine the optimal values of $v_i$'s, we do the following:
for given values of $n_i, c_i, c_i^*$, and $M_i$, $i=1, \ldots, L$,  and $c^{**}$,  we evaluate $Cost_1(\mathbf{v})$ for all possible choices of $v_1, \ldots, v_L$ such that for all $i$,  $n_iv_i\leq M_i$. Then the optimal values of $v_1, \ldots, v_L$ can be determined as the values for which the corresponding  mean  cost rate function $Cost_1(\mathbf{v})$ is minimum.
\begin{Remark}
  {\rm If the  system has a $k$-out-of-$n$ structure with independent components from multiple types components, then \eqref{1}  reduces to the result of \cite{Eryilmaz-main-2018}. This is so because for such system the survival signature  is obviously given as
\begin{align*}
\Phi(l_1,\ldots, l_L)=\left\{
\begin{tabular}{ll}
1 & $\sum_{j=1}^L l_j\geq k$\\
0 & otherwise,\\
\end{tabular}
\right.
\end{align*}
 i.e. the system works if at least $k$ components are alive.
}
\end{Remark}
\subsection{Cost function based on preventive  replacement}
\quad In this section, we propose a kind of age replacement preventive maintenance policy for the system with multiple types of components described in subsection 1.2.  The policy of renewing the system performed by the operator is such that it is replaced at failure time or at a predetermined time $\tau$, whichever occurs first.
There are many papers about age replacement strategy; the interested reader can see for example,  \cite{zhao2016}, \cite{ashrafi}    and \cite{nakagawa2019}.
\cite{Mannai2018}  found the optimal configuration of a $k$-out-of-$n$ system so that the
expected total costs of the system under some generalized age replacement policies are minimized.

\quad Here, suppose that the operator  has $M_i$ components of type $i$ as spare available, and he/she decides to add $v_i$ components to each of the components of type $i$, where $n_iv_i\leq M_i$.
Under the implemented policy, here, the aim is to find the optimal number of $v$'s such that the mean cost rate we impose below is minimized.

\quad  If the replacement occurs after the system failure, i.e., $T_R\leq \tau$, then, considering the costs  $c_i$,  $c_i^*$ and  $c^{**}$ as defined in the previous subsection, the average cost of renewing the system is obtained as
\begin{align*}
M_1(\mathbf{v})&=\sum_{i=1}^L c_i E(X_i(T_R)|T_R\leq \tau)+\sum_{i=1}^L c^*_i E(n_i(v_i+1)-X_i(T_R)|T_R\leq \tau)+c^{**}\nonumber\\
&=\sum_{i=1}^L (v_i+1)c_i E(X_i(T)|T\leq\tau)+\sum_{i=1}^L (v_i+1)c^*_i E(n_i-X_i(T)|T\leq\tau)+c^{**}\nonumber\\
&=\sum_{i=1}^L (c_i-c^*_i)(v_i+1) E(X_i(T)|T\leq\tau)+\sum_{i=1}^L (v_i+1)c^*_i n_i+c^{**},
\end{align*}
where $T$ is the lifetime of the system before redundancy allocation.

\quad If the system is replaced before failure, i.e., $T_R>\tau$, then, by the costs $c_i$ and $c_i^*$, $i=1, \ldots, L$  for renewing the failed components and refreshing  the alive components of type $i$, respectively,   the system will be  as good as the new condition.
Let $N_i(\tau)$ be the number of failed components of type $i$ on $[0, \tau]$. Then the average cost of renewing the system is defined as
\begin{align*}
M_2(\mathbf{v})&=\sum_{i=1}^L c_i E(N_i(\tau)|T_R>\tau)+\sum_{i=1}^L c^*_i E(n_i(v_i+1))-N_i(\tau)|T_R>\tau)\nonumber\\
&=\sum_{i=1}^L (c_i-c^*_i)(v_i+1) E(N_i(\tau)|T>\tau)+\sum_{i=1}^L (v_i+1)c^*_i n_i.
\end{align*}
Consequently, the mean cost rate function of the system renewing at time $\min(\tau, T_R)$ is achieved as
\begin{align}\label{cost2}
Cost_2(\mathbf{v})=\frac{M_1(\mathbf{v})P(T_R\leq \tau)+M_2(\mathbf{v})P(T_R>\tau)}{E(\min(\tau,T_R))},
\end{align}
{where it is attained that}
\begin{align*}
E(\min(\tau,T_R))=\int_{0}^{\tau}\bar{F}_{T_R}(y)dy.
\end{align*}
{For computing \eqref{cost2}}, we need to calculate  $E(N_i(\tau)|T>\tau)$ and $E(X_i(T)|T\leq \tau)$. For the first one, we have
\begin{align}
\label{EN12}
&E(N_i(\tau)|T>\tau)=\frac{1}{\bar{F}_T(\tau)}\sum_{j_i=0}^{n_i} j_i P(N_i(\tau)=j_i, T>\tau)\nonumber\\
&=\frac{1}{\bar{F}_T(\tau)}\sum_{j_1=0}^{n_1}\ldots\sum_{j_L=0}^{n_L} j_i P(T>\tau | N_1(\tau)=j_1,\ldots, N_L(\tau)=j_L)P(N_1(\tau)=j_1,\ldots, N_L(\tau)=j_L)\nonumber\\
&=\frac{1}{\bar{F}_T(\tau)}\sum_{j_1=0}^{n_1}\ldots\sum_{j_L=0}^{n_L} j_i \Phi(n_1-j_1,\ldots, n_L-j_L)\binom{n_1}{j_1}\ldots \binom{n_L}{j_L}B(\tau,j_1,\dots,j_L)
\end{align}
where
\begin{align}\label{Btau}
B(\tau,j_1,\dots,j_L)&=P(T_1^{(1)}\leq\tau, \ldots, T_{j_1}^{(1)}\leq\tau, T_{j_1+1}^{(1)}>\tau, \ldots, T_{n_1}^{(1)}>\tau, \nonumber\\
&\hspace{4cm}\ldots, T_1^{(L)}\leq\tau, \ldots, T_{j_L}^{(L)}\leq\tau, T_{j_L+1}^{(L)}>\tau, \ldots, T_{n_L}^{(L)}>\tau ).
\end{align}

Using a similar method as used in Lemma \ref{ui},   we  can calculate $E(X_i(T)|T\leq \tau)$, $i=1,\dots,L$, as follows.
\begin{align*}
E(X_i(T)|T\leq \tau)
&=n_iP(T_1^{(i)}\leq T|T\leq \tau)=n_i\frac{P(T_1^{(i)}\leq T,T\leq \tau)}{1-P(T>\tau)}, \quad j=1,\dots,L.
\end{align*}
Now, we can write
\begin{align*}
P(T_1^{(i)}\leq T,T\leq \tau)=\int_{0}^{\tau}\lim_{\delta\rightarrow 0}\frac{P(s<T\leq \tau, s<T_1^{(i)}\leq s+\delta)}{\delta}ds,\quad i=1,\dots,L,
\end{align*}
for which we have
\begin{small}
\begin{align*}
&P(s<T\leq \tau, s<T_1^{(i)}\leq s+\delta)\nonumber\\
&=\sum_{m_1=0}^{n_1}\ldots\sum_{m_i=0}^{n_i-1}\ldots\sum_{m_L=0}^{n_L}~\sum_{l_1=0}^{m_1}\ldots\sum_{l_L=0}^{m_L}P(s<T\leq \tau|s<T_1^{(i)}\leq s+\delta, C_1(\tau)=l_1,\ldots, C_i(\tau)=l_i\nonumber\\
&\hspace{1cm}C_L(\tau)=l_L, C_1(s)=m_1,\ldots, C_i(s)=m_i,\ldots, C_L(s)=m_L)P\big(s<T_1^{(i)}\leq s+\delta,C_1(\tau)=l_1,\nonumber\\
&\hspace{1cm},\ldots, C_i(\tau)=l_i,C_L(\tau)=l_L, C_1(s)=m_1,\ldots, C_i(s)=m_i,\ldots, C_L(s)=m_L)\big )\nonumber\\
&=\sum_{m_1=0}^{n_1}\ldots\sum_{m_i=0}^{n_i-1}\ldots\sum_{m_L=0}^{n_L}~\sum_{l_1=0}^{m_1}\ldots\sum_{l_L=0}^{m_L}[\Phi(m_1,\ldots, m_L)-\Phi(l_1,\ldots, l_L)]\left[\prod_{j=1,j\neq i}^L\binom{n_j}{m_j}\binom{m_j}{l_j}\right]\nonumber\\
&\times \binom{n_i-1}{m_i}\binom{m_i}{l_i}A_{\mathbf{m}, \mathbf{l}}^{(i)}(s,s+\delta, \tau),
\end{align*}
\end{small}
where
\begin{small}
\begin{align}\label{Asdelta}
&A_{\mathbf{m}, \mathbf{l}}^{(i)}(s,s+\delta, \tau)\nonumber\\
&=P(T_1^{(1)}>\tau, \ldots, T_{l_1}^{(1)}>\tau, s<T_{l_1+1}^{(1)}\leq \tau, \ldots, s<T_{m_1}^{(1)}\leq \tau, T_{m_1+1}^{(1)}\leq s,\ldots, T_{n_1}^{(1)}\leq s, \nonumber \\
&\hspace{0.2cm}\ldots, T_1^{(i)}>\tau, \ldots, T_{l_i}^{(i)}>\tau, s<T_{l_i+1}^{(i)}\leq\tau, \ldots, s<T_{m_i}^{(i)}\leq\tau, s<T_{m_i+1}^{(i)}\leq s+\delta, T_{m_i+2}^{(i)}\leq s,\ldots, T_{n_i}^{(i)}\leq s,
\nonumber \\
&\ldots,  T_1^{(L)}>\tau, \ldots, T_{l_L}^{(L)}>\tau, s<T_{l_L+1}^{(L)}\leq \tau, \ldots, s<T_{m_L}^{(L)}\leq \tau, T_{m_L+1}^{(L)}\leq s,\ldots, T_{n_L}^{(L)}\leq s ).
\end{align}
\end{small}
\vspace{0.1cm}
In the following theorem, the  probabilities in  \eqref{Btau} and \eqref{Asdelta} are represented based on the survival copula of components lifetimes.
\begin{Theorem}\label{Aml}
Using the inclusion-exclusion rule, we get the following expressions for $B(\tau,j_1,\dots,j_L)$ and $A_{\mathbf{m}, \mathbf{l}}^{(i)}(s,s+\delta, \tau)$, respectively.
\begin{align*}
B(\tau,j_1,\dots,j_L)&
={\sum_{b_1=0}^{j_1}\ldots\sum_{b_L=0}^{j_L}}(-1)^{b_1+\ldots+ b_L}\binom{j_1}{b_1}\ldots \binom{j_L}{b_L} \hat{C}(\underbrace{\bar{F}_1(\tau)}_{n_1-j_1+b_1},\underbrace{1}_{j_1-b_1},\ldots,\underbrace{\bar{F}_L(\tau)}_{n_L-j_L+b_L},\underbrace{1}_{j_L-b_L)})
\end{align*}
and
\begin{small}
\begin{align*}
&A_{\mathbf{m}, \mathbf{l}}^{(i)}(s,s+\delta, \tau)={\sum_{j_1=0}^{n_1-m_1}\ldots\sum_{j_i=0}^{n_i-m_i-1}\ldots\sum_{j_L=0}^{n_L-m_L}}(-1)^{j_1+\ldots + j_L}\binom{n_1-m_1}{j_1}\ldots \binom{n_i-m_i-1}{j_i}\ldots\binom{n_L-m_L}{j_L}\\
&\times {\sum_{d_1=0}^{m_1-l_1}\ldots\sum_{d_i=0}^{m_i-l_i}\ldots\sum_{d_L=0}^{m_L-l_L}}(-1)^{d_1+\ldots + d_L}\binom{m_1-l_1}{d_1}\ldots \binom{m_L-l_L}{d_L}\\
&\times \big[\hat{C}(\underbrace{\bar{F}_k(\tau)}_{l_k+d_k}, \underbrace{\bar{F}_k(s)}_{m_k-l_k+j_k-d_k},\underbrace{1}_{n_k-m_k-j_k}, for~ k=1,\ldots, L, k\neq i, \underbrace{\bar{F}_i(\tau)}_{l_i+d_i},  \underbrace{\bar{F}_i(s)}_{m_i-l_i+j_i-d_i+1}, \underbrace{1}_{n_i-m_i-j_i-1})\\
&-\hat{C}(\underbrace{\bar{F}_k(\tau)}_{l_k+d_k}, \underbrace{\bar{F}_k(s)}_{m_k-l_k+j_k-d_k},\underbrace{1}_{n_k-m_k-j_k}, for~ k=1,\ldots, L, k\neq i,\underbrace{\bar{F}_i(\tau)}_{l_i+d_i}, \underbrace{\bar{F}_i(s)}_{m_i-l_i+j_i-d_i}, \bar{F}_i(s+\delta), \underbrace{1}_{n_i-m_i-j_i-1} )\big].
\end{align*}
\end{small}
\end{Theorem}
\begin{proof}
See the Appendix
\end{proof}
\begin{Corollary}
For the particular case of  independent components, it can be deduced that
\begin{align}\label{2}
E(N_i(\tau)|T>\tau)=\frac{1}{\bar{F}(\tau)}\sum_{j_1=0}^{n_1}\ldots\sum_{j_L=0}^{n_L} j_i\Phi(n_1-j_1,\ldots, n_L-j_L)\prod_{l=1}^L \binom{n_l}{j_l}F_l^{j_l}(\tau)\bar{F}_l^{n_l-j_l}(\tau).
\end{align}
Also, in this case, we have
\begin{align*}
A_{\mathbf{m}, \mathbf{l}}^{(i)}(s,s+\delta, \tau)&=\big\{\prod_{j=1,j\neq i}^L F_j^{n_j-m_j}(s)[\bar{F}_j(\tau)-\bar{F}_j(s)]^{m_j-l_j}\bar{F}_j^{l_j}(\tau)\big\}\\
&\times[\bar{F}_i(s)-\bar{F}_i(s+\delta)]F_i^{n_i-m_i-1}(s)[\bar{F}_i(\tau)-\bar{F}_i(s)]^{m_i-l_i}\bar{F}_i^{l_i}(\tau)
\end{align*}
which, in turn, implies that
\begin{small}
\begin{align}\label{22}
&E(X_i(T)|T\leq \tau)=\frac{n_i}{1-\bar{F}_{T}(\tau)}\sum_{m_1=0}^{n_1}\ldots\sum_{m_i=0}^{n_i-1}\ldots\sum_{m_L=0}^{n_L}
\sum_{l_1=0}^{m_1}\ldots\sum_{l_L=0}^{m_L}
[\Phi(m_1,\ldots, m_L)-\Phi(l_1,\ldots, l_L)]\nonumber\\
&\times\left[\prod_{j=1,j\neq i}^L\binom{n_j}{m_j}\binom{m_j}{l_j}\right] \binom{n_i-1}{m_i}\binom{m_i}{l_i}\int_{0}^{\tau}\big\{\prod_{j=1,j\neq i}^L F_j^{n_j-m_j}(s)[\bar{F}_j(\tau)-\bar{F}_j(s)]^{m_j-l_j}\bar{F}_j^{l_j}(\tau)\big\}\nonumber\\
&\hspace{7.5cm}\times F_i^{n_i-m_i-1}(s)[\bar{F}_i(\tau)-\bar{F}_i(s)]^{m_i-l_i}\bar{F}_i^{l_i}(\tau)dF_i(s).
\end{align}
\end{small}
\end{Corollary}
\quad It is worth noting that if the structure of the system is $k$-out-of-$n$  whose components are  independent, then \eqref{2} and \eqref{22} are reduced to the results appeared in \cite{Eryilmaz-main-2018}.

\quad For given values $n_i, c_i, c_i^*, M_i$, $i=1,\ldots, L$,  $c^{**}$ and $\tau$, we aim to determine the optimal values of  $v_i$'s under the constraints $n_iv_i\leq M_i,$ $ i=1, \ldots, L$, such that the mean cost rate function  $Cost_2$ is minimized.

In the following, we present two examples to examine the aforementioned theoretical results.

\begin{Example}\label{eg11}
{\rm Consider the system depicted in Figure \ref{fig1}, given in \cite{feng2016} and \cite{eriylmaz-dependent-2018}. The system consists of six components in which components 1, 2, and 5 are of type one, and components 3, 4, and 6 are of type two. The survival signature of the system is presented in Table \ref{table1}.
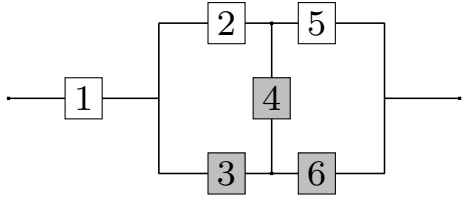
\begin{figure}[h!]
\centerline{\begin{tikzpicture}[node distance = 5 cm]
  \tikzset{LabelStyle/.style =   {scale=2, fill= white, text=black}}
     \node[shape = rectangle,draw, fill= black, text= black, inner sep =2 pt, outer sep= 0 pt, scale=0.2](A) at (-2,0) {};
     \node[shape = rectangle,draw, fill= white, text= black, inner sep =2 pt, outer sep= 0 pt, scale=1.7](B) at (-1,0) {\scriptsize 1};
     \node[shape = rectangle,draw, fill= lightgray, text= black, inner sep =2 pt, outer sep= 0 pt, scale=0.01](D) at (0,0) {};
     \node[shape = rectangle,draw, fill= lightgray, text= black, inner sep =2 pt, outer sep= 0 pt, scale=0.01](E) at (0,1) {};
     \node[shape = rectangle,draw, fill= lightgray, text= black, inner sep =2 pt, outer sep= 0 pt, scale=0.01](F) at (0,-1) {};
      \node[shape = rectangle,draw, fill= white, text= black, inner sep =2 pt, outer sep= 0 pt, scale=1.7](G) at (0.9,1) {\scriptsize 2};
      \node[shape = rectangle,draw, fill= black, text= black, inner sep =2 pt, outer sep= 0 pt, scale=0.2](GH) at (1.5,1) {};
      \node[shape = rectangle,draw, fill= white, text= black, inner sep =2 pt, outer sep= 0 pt, scale=1.7](H) at (2.1,1) {\scriptsize 5};
       \node[shape = rectangle,draw, fill= lightgray, text= black, inner sep =2 pt, outer sep= 0 pt, scale=0.01](I) at (3,1) {};
        \node[shape = rectangle,draw, fill= lightgray, text= black, inner sep =2 pt, outer sep= 0 pt, scale=1.7](DD) at (1.5,0) {\scriptsize 4};
        \node[shape = rectangle,draw, fill= lightgray, text= black, inner sep =2 pt, outer sep= 0 pt, scale=1.7](J) at (0.9,-1) {\scriptsize 3};
     \node[shape = rectangle,draw, fill= black, text= black, inner sep =2 pt, outer sep= 0 pt, scale=0.2](JK) at (1.5,-1) {};
    \node[shape = rectangle,draw, fill= lightgray, text= black, inner sep =2 pt, outer sep= 0 pt, scale=1.7](K) at (2.1,-1) {\scriptsize 6};
     \node[shape = rectangle,draw, fill= lightgray, text= black, inner sep =2 pt, outer sep= 0 pt, scale=0.01](II) at (3,-1) {};
         \node[shape = rectangle,draw, fill= lightgray, text= black, inner sep =2 pt, outer sep= 0 pt, scale=0.01](L) at (3,0) {};
       \node[shape = rectangle,draw, fill= black, text= black, inner sep =2 pt, outer sep= 0 pt, scale=0.2](M) at (4,0) {};
     \draw[semithick](A) -- (B);
     \draw[semithick](B) -- (D);
     \draw[semithick](D) -- (E);
     \draw[semithick](D) -- (F);
      \draw[semithick](E) -- (G);
      \draw[semithick](GH) -- (DD);
     \draw[semithick](JK) -- (DD);
     \draw[semithick](G) -- (H);
     \draw[semithick](H) -- (I);
     \draw[semithick](I) -- (L);
     \draw[semithick](F) -- (J);
     \draw[semithick](J) -- (K);
     \draw[semithick](II) -- (L);
      \draw[semithick](K) -- (II);
     \draw[semithick](L) -- (M);
  \end{tikzpicture}}
   \caption{\small  The system in Example \ref{eg11} with two types of components}
   \hfill
   \label{fig1}
\end{figure}

\begin{table}[h]
\small
\caption{Survival signature of the system in Figure \ref{fig1}}\label{table1}
\begin{center}
\begin{tabular}{cccccc}
\hline
$l_1$ & $l_2$ & $\Phi(l_1, l_2)$ & $l_1$ & $l_2$ & $\Phi(l_1, l_2)$\\
\hline
0 & 0 & 0 & 2 & 0 & 0\\
0 & 1 & 0 & 2 & 1 & 0\\
0 & 2 & 0 & 2 & 2 & $\frac{4}{9}$\\
0 & 3 & 0 & 2 & 3 & $\frac{2}{3}$\\
1 & 0 & 0 & 3 & 0 & 1\\
1 & 1 & 0 & 3 & 1 & 1\\
1 & 2 & $\frac{1}{9}$ & 3 & 2 & 1\\
1 & 3 & $\frac{1}{3}$ & 3 & 3 &1\\
\hline
\end{tabular}
\end{center}
\end{table}
Assume that the dependency structure of the component lifetimes  is modeled by  a parametric family of copulas known as  Gumbel-Hougaard family defined as
 \begin{align*}
&\hat{C}(u_1,\ldots, u_n)=\exp\left(-\left[(-\ln u_1)^{\alpha})+\ldots+(-\ln u_n)^{\alpha}) \right] ^{1/\alpha} \right),
 \end{align*}
 where $\alpha\geq 1$ is the dependency parameter in the family. The value $\alpha=1$ corresponds to the independent condition.  Let the component lifetimes of the two types follow exponential distributions with reliability functions $\bar{F}_i(t)=e^{-t\theta_i}$, where we assume that $\theta_1=0.2$,  and $\theta_2=0.3$. If there are  $M_1=9$ and $M_2=6$ components from type 1 and type 2 as spares, respectively, then $v_1\in\{0,1,2,3\}$ and $ v_2\in\{0,1,2\}$.
To find the optimal number of redundant  components for each type, we use the following values for the replacement costs: $c_1=3, c_2=2, c_1^*=1.5, c_2^*=1$, and $ c^{**}=10$.
{
For computing the numerator of \eqref{poi}, we need to compute $E(X_i(T)), i=1, 2$.  From Lemma \ref{ui} and Theorem \ref{Am}, we have
\begin{small}
\begin{align*}
&E(X_1(T))=n_1\int_{0}^{\infty} \lim_{\delta\rightarrow 0}\frac{1}{\delta}\sum_{m_1=0}^{n_1-1}\sum_{m_2=0}^{n_2}\Phi(m_1+1, m_2)  \binom{n_1-1}{m_1}\binom{n_2}{m_2}A_{\mathbf{m}}^{(1)}(t,\delta)dt,
\end{align*}
\end{small}
where
\begin{align*}
&A_{\mathbf{m}}^{(1)}(t, \delta)={\sum_{j_1=0}^{n_1-m_1-1} \sum_{j_2=0}^{n_2-m_2}(-1)^{j_1+ j_2}\binom{n_1-m_1-1}{j_1}\binom{n_2-m_2}{j_2}} \\
&\times  \Big[ e^{-[(m_1+j_1+1)(t \theta_1)^{\alpha}+(m_2+j_2)(t \theta_2)^{\alpha}]^{1/\alpha}}-e^{-[(m_1+j_1)(t \theta_1)^{\alpha}+((t+\delta)\theta_1)^{\alpha}+(m_2+j_2)(t \theta_2)^{\alpha}]^{1/\alpha}}\Big].
\end{align*}
Thus, we get
\begin{align*}
&E(X_1(T))\\
&=n_1 \sum_{m_1=0}^{n_1-1}\sum_{m_2=0}^{n_2}\binom{n_1-1}{m_1}\binom{n_2}{m_2}\Phi(m_1+1, m_2)\sum_{j_1=0}^{n_1-m_1-1} \sum_{j_2=0}^{n_2-m_2} (-1)^{j_1+ j_2} 
\binom{n_1-m_1-1}{j_1}\binom{n_2-m_2}{j_2}\\
&\times \theta_1^{\alpha}\big[(m_1+j_1+1)\theta_1^{\alpha}+(m_2+j_2)\theta_2^{\alpha}\big]^{-1}.
\end{align*}
Similarly, we have
\begin{align*}
&E(X_2(T))\\
&=n_2 \sum_{m_1=0}^{n_1}\sum_{m_2=0}^{n_2-1}\binom{n_1}{m_1}\binom{n_2-1}{m_2}\Phi(m_1, m_2+1)\sum_{j_1=0}^{n_1-m_1} \sum_{j_2=0}^{n_2-m_2-1} (-1)^{j_1+ j_2}
\binom{n_1-m_1}{j_1}\binom{n_2-m_2-1}{j_2}\\
&\times\theta_2^{\alpha}\big[(m_1+j_1)\theta_1^{\alpha}+(m_2+j_2+1)\theta_2^{\alpha}\big]^{-1}.
\end{align*}
Also,  the denominator of \eqref{poi} can be written as follows 
\begin{align*}
E(T_R)&=\int_0^{\infty} \bar{F}_{T_R}(t)  dt={\sum_{l_1=0}^{n_1}\sum_{l_2=0}^{n_2}}\binom{n_1}{l_1} \binom{n_2}{l_2}\Phi(l_1, l_2){\sum_{i_1=0}^{n_1-l_1}\sum_{i_2=0}^{n_2-l_2}}(-1)^{i_1+ i_2}\binom{n_1-l_1}{i_1} \binom{n_2-l_2}{i_2} \nonumber\\
&\times
\int_0^{\infty} e^{-\big[(i_1+l_1)\big(-\ln(1-(1-\exp[-t \theta_1])^{v_1+1})\big)^{\alpha}+(i_2+l_2)\big(-\ln(1-(1-\exp[-t \theta_2])^{v_2+1})\big)^{\alpha} \big]^{1/\alpha}} dt,
\end{align*}
which should be evaluated numerically by  suitable softwares such as $Mathematica$.
}
The values of $Cost_1(\mathbf{v})$ for different combinations of $v_1$ and $ v_2$ are presented in Table \ref{table11} for two values  $\alpha=2$ (dependent components) and $\alpha=1$ (independent components). { It is seen that  $v_1=2$ and $v_2=0$ are the optimal choices for the number of  redundant components of the first and the second types,  respectively, under the criterion  $Cost_1$ in the case $\alpha=2$, and $v_1=3$ and $v_2=0$  are the optimal numbers in the case $\alpha=1$.  }
\begin{table}[h!]
\small
\begin{center}
\caption{{The values of $Cost_1(\mathbf{v})$  and $Cost_2(\mathbf{v})$   in Example \ref{eg11} in the case of dependent components $(\alpha=2)$ and independent components  $(\alpha=1)$.}}\label{table11}
\end{center}
\begin{center}
\begin{tabular}{|cc|c|c|c|c|}
\hline
$v_1$ & $v_2$ & {$Cost_1(v_1, v_2)$}  & {$Cost_2(v_1, v_2)$} & {$Cost_1(v_1, v_2)$}  & { $Cost_2(v_1, v_2)$}\\
&&$\alpha=2$& $\alpha=2$&$\alpha=1$& $\alpha=1$
\\
\hline
0 & 0 & 6.33927 & 8.2455 & 9.36071 & 9.70214\\
0 & 1 & 6.81922 & 9.44774 & 9.38725 & 10.3790\\
0 & 2 & 7.5289 & 11.3022 & 9.87544 & 11.9363\\
1 & 0 & 5.20719 & {\bf 7.77258} & 7.09069 & {\bf 8.48716}\\
1 & 1 & 5.82298 & 9.2981 & 7.51217 & 9.7885\\
1 & 2 & 6.35331 & 10.7302 & 7.98448 & 11.3693\\
2 & 0 & {\bf 4.99041} & 9.13115 & 6.56325 & 9.88167\\
2 & 1 & 5.58924 & 10.7299 & 7.06002 & 11.4307\\
2 & 2 & 6.13251 & 12.28921 & 7.53562 & 13.0756\\
3 & 0 & 5.04518 & 11.137 & {\bf 6.48197} & 12.035\\
3 & 1 & 5.59315 & 12.7178 & 6.98476 & 13.6766\\
3 & 2 & 6.11375 & 14.2941 & 7.45447 & 15.3508\\
\hline
\end{tabular}
\end{center}
\end{table}

Suppose that the described system is maintained under the aforementioned age replacement policy, where we assume that  $\tau=2$, i.e., the  replacement time of the system is $\min(T_R, 2)$.
{
From equations \eqref{EN12} and \eqref{Btau} the  mean number of failed components of $i$th type at time $\tau$, before system failure, is evaluated by the following expression. 
\begin{align*}
&E(N_i(\tau)|T>\tau)=\frac{1}{\bar{F}_T(\tau)}\sum_{j_1=0}^{n_1}\sum_{j_2=0}^{n_2} j_i \Phi(n_1-j_1, n_2-j_2)\binom{n_1}{j_1} \binom{n_2}{j_2}
{\sum_{b_1=0}^{j_1}\sum_{b_2=0}^{j_2}}(-1)^{b_1+ b_2}\binom{j_1}{b_1} \binom{j_2}{b_2} \\
&\times \exp[-\tau\big( \theta_1^{\alpha}(n_1-j_1+b_1)+\theta_2^{\alpha}(n_2-j_2+b_2) \big)^{1/\alpha}], ~i=1, 2, 
\end{align*}
where from \eqref{Fbar-dep}, we get
 \begin{align*}
\bar{F}_T(\tau)&={\sum_{l_1=0}^{n_1}\sum_{l_2=0}^{n_2}}{\sum_{i_1=0}^{n_1-l_1}\sum_{i_2=0}^{n_2-l_2}}(-1)^{i_1+ i_1}\binom{n_1}{l_1} \binom{n_2}{l_2}\binom{n_1-l_1}{i_1} \binom{n_2-l_2}{i_2}\Phi(l_1, l_2)\nonumber\\
&\times \exp[-\tau\big((i_1+l_1)\theta_1^{\alpha}+(i_2+l_2)\theta_2^{\alpha} \big)^{1/\alpha}].
\end{align*}
Next, for the mean number of failed components at the time of the system failure given that the system has failed before $\tau$, we have
\begin{small}
\begin{align*}
&E(X_1(T)|T\leq \tau)\\
&=\frac{n_1}{1-\bar{F}_T(\tau)}\sum_{m_1=0}^{n_1-1}\sum_{m_2=0}^{n_2}\sum_{l_1=0}^{m_1}\sum_{l_2=0}^{m_2}[\Phi(m_1, m_2)-\Phi(l_1, l_2)]\binom{n_1-1}{m_1}\binom{m_1}{l_1}\binom{n_2}{m_2}\binom{m_2}{l_2}\\
&\times 
{\sum_{j_1=0}^{n_1-m_1-1}\sum_{j_2=0}^{n_2-m_2}}(-1)^{j_1+ j_2}\binom{n_1-m_1-1}{j_1} \binom{n_2-m_2}{j_2} {\sum_{d_1=0}^{m_1-l_1}\sum_{d_2=0}^{m_2-l_2}}(-1)^{d_1+d_2}\binom{m_1-l_1}{d_1} \binom{m_2-l_2}{d_2}\\
&\times \frac{\theta_1^{\alpha}\big[ e^{-\tau\big((l_1+d_1)\theta_1^{\alpha}+(l_2+d_2)\theta_2^{\alpha}  \big)^{1/\alpha}}- 
 e^{-\tau\big((m_1+j_1+1)\theta_1^{\alpha}+(m_2+j_2)\theta_2^{\alpha}  \big)^{1/\alpha}}
\big]}{(m_1-l_1+j_1-d_1+1)\theta_1^{\alpha}+(m_2-l_2+j_2-d_2)\theta_2^{\alpha}},
\end{align*}
\end{small}
and
\begin{small}
\begin{align*}
&E(X_2(T)|T\leq \tau)\\
&=\frac{n_2}{1-\bar{F}_T(\tau)}\sum_{m_1=0}^{n_1}\sum_{m_2=0}^{n_2-1}\sum_{l_1=0}^{m_1}\sum_{l_2=0}^{m_2}[\Phi(m_1, m_2)-\Phi(l_1, l_2)]\binom{n_1}{m_1}\binom{m_1}{l_1}\binom{n_2-1}{m_2}\binom{m_2}{l_2}\\
&\times 
{\sum_{j_1=0}^{n_1-m_1}\sum_{j_2=0}^{n_2-m_2-1}}(-1)^{j_1+ j_2}\binom{n_1-m_1}{j_1} \binom{n_2-m_2-1}{j_2} {\sum_{d_1=0}^{m_1-l_1}\sum_{d_2=0}^{m_2-l_2}}(-1)^{d_1+d_2}\binom{m_1-l_1}{d_1} \binom{m_2-l_2}{d_2}\\
&\times \frac{\theta_2^{\alpha}\big[ e^{-\tau\big((l_1+d_1)\theta_1^{\alpha}+(l_2+d_2)\theta_2^{\alpha}  \big)^{1/\alpha}}- 
 e^{-\tau\big((m_1+j_1)\theta_1^{\alpha}+(m_2+j_2+1)\theta_2^{\alpha}  \big)^{1/\alpha}}
\big]}{(m_1-l_1+j_1-d_1)\theta_1^{\alpha}+(m_2-l_2+j_2-d_2+1)\theta_2^{\alpha}}.
\end{align*}
\end{small}
By substituting these results in \eqref{cost2}, the mean cost rate of replacement strategy can be evaluated. 
}
 \quad The values of  $Cost_2(\mathbf{v})$ are calculated  for different combinations of $v_1$ and $v_2$, in Table \ref{table11}  for both dependent and independent situations. It follows, from the results of the table, that for $v_1=1$ and $v_2=0$, the mean cost rate $Cost_2(\mathbf{v})$ is minimized, in both cases $\alpha=1,2$.
 
 {\quad  In order to investigate the robustness of our strategies  concerning the model parameters, we calculate some numerical results based on these parameters.}
The results in Table \ref{table66} shows the effect of the dependency parameter $\alpha$ on
 the optimal values of $v_1$ and $v_2$,  for different values of $\alpha$.  As seen, when $\alpha$ increases (i.e., we get far from independence)  the number of redundant components decreases { based on the objective function $Cost_1(v_1,v_2)$, but remain unchanged under the $Cost_2(v_1,v_2)$. }{This makes sense since under the more dependency the {$MTTF$} is increased, and the need to spare components reduces.
Also it should be noted that the higher the $\alpha$, the lower mean cost rates. }
\begin{table}[h!]
\small
\caption{\small The optimum values of  $\mathbf{v}$ by minimizing $Cost_i(v_1,v_2)$, $i=1,2$, for different $\alpha$ in Example \ref{eg11}}\label{table66}
\begin{center}
{
\begin{tabular}{|c|ccccccccccc|}
\hline
$\alpha$& 1 & 1.2 & 1.4 & 1.6 & 1.8 & 2 & 2.2 & 2.4 & 2.6 & 2.8 & 3 \\ \hline
$v_1$ & 3 & 3 & 2 & 2 & 2 & 2 & 2 & 2 & 2 & 2 & 2  \\
$v_2$ & 0 & 0 & 0 & 0 & 0 & 0 & 0 & 0 & 0 & 0 & 0  \\
$Cost_1(v_1,v_2)$ & 6.48 & 5.98 & 5.63 & 5.36 & 5.15  & 4.99 & 4.86 & 4.75 &  4.66 & 4.58 & 4.52\\
 \hline
$v_1$ & 1 & 1 & 1 & 1 & 1 & 1 & 1 & 1 & 1 & 1 & 1  \\
$v_2$ & 0 & 0 & 0 & 0 & 0 & 0 & 0 & 0 & 0 & 0 & 0  \\
${Cost_2}(v_1,v_2)$ & 8.48  & 8.28 & 8.11 & 7.98 & 7.87 & 7.77 & 7.69 & 7.62 & 7.57 & 7.50 & 7.46\\
\hline
\end{tabular}
}
\end{center}
\end{table}

 \quad For exploring the sensitivity  of  the proposed models with respect to components costs,
 $\mathbf{c}=(c_1,c_2)$ and $\mathbf{c}^*=(c_1^*,c_2^*)$, we have provided some numerical results in Table \ref{table412}, for $\alpha=2$.   In the top part of the left panel of the table, we observe that for fixed values of $c_1^*=1.5$ and $c_2^*=1$, the increase in the costs $c_1$ and $c_2$ results a reduction to the number of optimal values of $v_1$ and $v_2$.  In the bottom  part of the left panel of the table, the costs $c_i$ and $c^*_i$, $i=1,2$, of the two types are swapped. In this case,  when the costs $c_1^*$ and $c_2^*$ are fixed as $c_1^*=1$ and $c_2^*=1.5$, then we again observe that the increase in the costs $c_1$ and $c_2$ results a decline to the number of optimal values of $v_1$ and $v_2$. In the top part of the right panel of Table \ref{table412}, it can be seen that for fixed values of renewing failed components as $c_1=6, c_2=5.5$,  the decrease in the costs $c_1^*$ and $c_2^*$ results an increase to the number of optimal values of $v_1$ and $v_2$. As shown in the bottom part of the right panel, the same result holds by swapping the costs of the components of type 1 and type 2. { It is seen that in all four parts of Table \ref{table412}
the value of costs and the numbers of spare components are inversely related to each other.}
\begin{table}[h!]
\small
\caption{The optimum values of  $\mathbf{v}$ by minimizing $Cost_1$ for different costs  in Example \ref{eg11}}\label{table412}
\begin{center}
\begin{tabular}{|ccccc||ccccc|}
\hline
$\mathbf{c}$ & $\mathbf{c^*}$ & { $v_1$} & {$v_2$} & { $Cost_1$} & $\mathbf{c}$ & $\mathbf{c^*}$ & { $v_1$} & {$v_2$} & {$ Cost_1$}\\
\hline                                           
(1.6, 1.1) & (1.5,1) & 3 & 0 & 4.0492 & (6, 5.5) & (5.9, 5.4) & 1 & 0 & 11.9478\\
(1.7, 1.2) & (1.5,1) & 3 & 0 & 4.1278 & (6, 5.5) & (5.7, 5.2) & 1 & 0 & 11.7509\\
(1.8, 1.3) & (1.5,1) & 3 & 0 & 4.2065 & (6, 5.5) & (5.5, 5.0) & 1 & 0 & 11.5540\\
(1.9, 1.4) & (1.5,1) & 2 & 0 & 4.2845 & (6, 5.5) & (5.3, 4.8) & 1 & 0 & 11.3571\\
(2, 1.5) & (1.5,1) & 2 & 0 & 4.3599 & (6, 5.5) &   (5.1, 4.6) & 2 & 0 & 11.1543\\
(2.5, 2) & (1.5,1) & 2 & 0 & 4.7369 & (6, 5.5) &  (5, 4.5) & 2 & 0 & 11.0493 \\
(3, 2.5) & (1.5,1) & 2 & 0 & 5.1139 & (6, 5.5) &  (4.5, 4) & 2 & 0 & 10.5245\\
(3.5, 3) & (1.5,1) & 2 & 0 & 5.4908 & (6, 5.5) &  (4, 3.5) & 2 & 0 & 9.9997\\
(4, 3.5) & (1.5,1) & 2 & 0 & 5.8678 & (6, 5.5) &  (3.5, 3) & 2 & 0 & 9.4750 \\
(4.5, 4) & (1.5,1) & 2 & 0 & 6.2448 & (6, 5.5) &  (3, 2.5) & 2 & 0 & 8.9502\\
\hline
(1.5, 2) & (1, 1.5) & 3 & 0 & 3.7874 & (5.5, 6) & (5, 5.5) & 2 & 0 & 11.1232\\
(2, 2.5) & (1, 1.5) & 3 & 0 & 4.1807 & (5.5, 6) & (4.5, 5) & 2 & 0 & 10.5984\\
(2.5, 3) & (1, 1.5) & 3 & 0 & 4.5740 & (5.5, 6) & (4, 4.5) & 2 & 0 & 10.0737\\
(3, 3.5) & (1, 1.5) & 3 & 0 & 4.9673 & (5.5, 6) & (3.5, 4) & 2 & 0 & 9.5489\\
(3.5, 4) & (1, 1.5) & 3 & 0 & 5.3606 & (5.5, 6) & (3, 3.5) & 2 & 0 & 9.0241\\
(4, 4.5) & (1, 1.5) & 3 & 0 & 5.7539 & (5.5, 6) & (2.5, 3) & 2 & 0 & 8.4992\\
(4.5, 5) & (1, 1.5) & 3 & 0 & 6.1472 & (5.5, 6) & (2, 2.5) & 2 & 0 & 7.9745\\
(5, 5.5) & (1, 1.5) & 3 & 0 & 6.5405 & (5.5, 6) & (1.5, 2) & 2 & 0 & 7.4497\\
(5.5, 6) & (1, 1.5) & 2 & 0 & 6.9249 & (5.5, 6) & (1, 1.5) & 2 & 0 & 6.9249\\
(6, 6.5) & (1, 1.5) & 2 & 0 & 7.0550 & (5.5, 6) & (0.5, 1) & 3 & 0 & 6.3664\\
(6.5, 7) & (1, 1.5) & 2 & 0 & 7.4320 & (5.5, 6) & (0, 0.5) & 3 & 0 & 5.7991\\
\hline
\end{tabular}
\end{center}
\end{table}
Table \ref{table44} shows the behavior of the number of redundant components from another point of view. In the left panel of the table, we have kept $c_2 (c^*_2)$ constant and have increased the values of $c_1 (c^*_1)$. In fact, we have assumed that $c_1=\omega c_2$ and $c_1^*=\omega c_2^*$ for $\omega=1.5,2,3, ..., 10$. As seen, when $\omega$ increases the optimal value of the redundant component $v_1$ decreases { and the optimal value of $v_2$ increases}. In the right panel of Table \ref{table44}, we exchange the costs of type 1 and 2, i.e. we assume $c_2=\omega c_1$ and $c_2^*=\omega c_1^*$. {In this case, we observe no changes in the number of redundant components $v_1, v_2$ when $\omega$ increases.}\\
\begin{table}[h!]
\small
\caption{The optimum values of  $\mathbf{v}$ by minimizing $Cost_1$ for different costs  in Example \ref{eg11}}\label{table44}
\begin{center}
\begin{tabular}{|ccccc||ccccc|}
\hline
$\mathbf{c}$ & $\mathbf{c^*}$ &{ $v_1$ }& {$v_2$ }& {$Cost_1$} &$\mathbf{c}$ & $\mathbf{c^*}$ & {$v_1$} & {$v_2$} &{ $Cost_1$} \\
\hline
(3,2) & (1.5,1) & 2 & 0 & 4.9904 &  (2,3) & (1,1.5) & 3 & 0 & 4.2859\\
(4,2) & (2,1) & 2 & 0 & 5.9203 & (2,4) & (1,2) & 3 & 0 & 4.5832\\
(6,2) & (3,1) & 1 & 0 & 7.5921 & (2,6) & (1,3) & 3 & 0 & 5.1778\\
(8,2) & (4,1) & 1 & 0 & 9.1824 & (2,8) & (1,4) & 3 & 0 & 5.7725\\
(10,2) & (5,1) & 0 & 1 & 10.6840 & (2,10) & (1,5) & 3 & 0 & 6.3671\\
(12,2) & (6,1) & 0 & 1 & 11.7883 & (2,12) & (1,6) & 3 & 0 & 6.9617\\
(14,2) & (7,1) & 0 & 1 & 12.8925 & (2,14) & (1,7) & 3 & 0 & 7.5563\\
(16,2) & (8,1) & 0 & 1 & 13.9967 & (2,16) & (1,8) & 3 & 0 & 8.1509\\
(18,2) & (9,1) & 0 & 1 & 15.1009 & (2,18) & (1,9) & 3 & 0 & 8.7456\\
(20,2) & (10,1) & 0 & 1 & 16.2052 & (2,20) & (1,10) & 3 & 0 & 9.3402\\
\hline
\end{tabular}
\end{center}
\end{table}

\quad {In this example, the distributions of the components lifetimes are ordered such that $\bar{F}_1(t)\geq \bar{F}_2(t)$, for all $t>0$, i.e. the reliability (and subsequently, the $MTTF$) of the components of type one is more than type two.
Note that, according to the system structure, it is reveal that the components of type one are generally in more critical positions than those of type two. 
Hence, one should intuitively expect that the optimal solution, according to the cost criterion, would be the case in which one allocates more components of type one than type two.}

\quad { To see whether this fact affects the number of $v_i$'s, we let $\bar{F}_2^*(t)=e^{-0.07 t^{2}}$, $t>0$,  be the reliability function of the components of type two. In this case, the two reliability functions cross each other such that $\bar{F}_1(t)<\bar{F}_2^*(t)$ for $t<2.86$ and $\bar{F}_1(t)>\bar{F}_2^*(t)$ for $t>2.86$. Note that the $MTTF$ for components of type two in this new case is the same as the previous one.
 By fixing the other parameters as before, we get the results given in Table \ref{table660}.
We see that although the distributions cross each other,  the optimal numbers of components in Table \ref{table660} are mostly the same as those in Table \ref{table66}, perhaps since the $MTTF$s have not been changed in both cases.
\begin{table}[h!]
\small
\caption{\small {The optimum values of  $\mathbf{v}$ by minimizing $Cost_1(v_1,v_2)$, for different values of $\alpha$ in Example \ref{eg11}}}\label{table660}
\begin{center}
{
\begin{tabular}{|c|ccccccccccc|}
\hline
$\alpha$& 1 & 1.2 & 1.4 & 1.6 & 1.8 & 2 & 2.2 & 2.4 & 2.6 & 2.8 & 3 \\ \hline
$v_1$ & 3 & 2 & 2 & 2 & 2 & 2 & 2 & 2 & 2 & 2 & 2  \\
$v_2$ & 0 & 0 & 0 & 0 & 0 & 0 & 0 & 0 & 0 & 0 & 0  \\
$Cost_1(v_1,v_2)$ & 6.53 & 6.06 & 5.71 & 5.45 & 5.25  & 5.08 & 4.95 & 4.85 &  4.76 & 4.68 & 4.62\\
 \hline
\end{tabular}
}
\end{center}
\end{table}
}

{
As a final point, to see the effect of survival copula on the optimal numbers of $v_1$ and $v_2$,  we suppose that the dependence structure is followed by the Clayton copula with the following form  
\begin{align*}
\hat{C}(u_1,\ldots, u_n)=\left(u_1^{-1/\alpha}+\ldots + u_n^{-1/\alpha}-n+1  \right)^{-\alpha}, \alpha>0.
\end{align*}
The parameter $\alpha$ manages the dependency degree of the copula, and the limiting case $\alpha=0$, gives the independence. We obtain the optimal values of redundant components for some values of $\alpha$ in Table \ref{table6687}. As it can be seen, by increasing $\alpha$, the values of $v_1$ and $v_2$ and also the mean cost rate show increment, which is in contradiction with the results in Table \ref{table66}. Hence, the output of the optimization problem strongly pertains to the functional structure of dependence, not only to the dependence parameter.}
\begin{table}[h!]
\small
\caption{\small {The optimum values of  $\mathbf{v}$ by minimizing $Cost_1(v_1,v_2)$ for different $\alpha$ under the Clayton copula  in Example \ref{eg11}}}\label{table6687}
\begin{center}
{
\begin{tabular}{|c|cccccc|}
\hline
$\alpha$ & 0.001 & 0.1 & 1 & 2 & 3  & 4   \\ \hline
$v_1$ & 2 & 2 & 2 & 2 & 3 & 3    \\
$v_2$ & 0 & 0 & 0 & 0 & 0 & 0    \\
$Cost_1(v_1,v_2)$ & 3.71 & 3.97 & 5.26 & 5.75 & 5.98  & 6.09  \\
\hline
\end{tabular}
}
\end{center}
\end{table}
}
\end{Example}
In the following example we consider an 8-component system consists of three types of components. For fixed values of $v_i$'s, we minimize the  function $Cost_1$ and also the function $Cost_2$ in the case that replacement time of unfailed system, $\tau$, is considered as the variable of interest.
\begin{Example}\label{eg2}
{\rm Consider the system depicted in Figure \ref{fig2}, given in \cite{heu}. The system has eight components from which three components (1, 2, and 3) are of type one, three components (4, 5, and 7) are of type two, and two components (6 and 8) are of type three. The values of the system survival signature are computed in \cite{heu}, and hence we refer the reader to the cited paper for the details.
\begin{figure}[h!]
\centerline{\begin{tikzpicture}[node distance = 5 cm]
  \tikzset{LabelStyle/.style =   {scale=2, fill= white, text=black}}
     \node[shape = rectangle,draw, fill= lightgray, text= black, inner sep =2 pt, outer sep= 0 pt, scale=0.01](D) at (0,0) {};
     \node[shape = rectangle,draw, fill= lightgray, text= black, inner sep =2 pt, outer sep= 0 pt, scale=0.01](E) at (0,1) {};
     \node[shape = rectangle,draw, fill= lightgray, text= black, inner sep =2 pt, outer sep= 0 pt, scale=0.01](F) at (0,-1) {};
      \node[shape = rectangle,draw, fill= white, text= black, inner sep =2 pt, outer sep= 0 pt, scale=1.7](G) at (0.9,1) {\scriptsize 1};
      \node[shape = rectangle,draw, fill= black, text= black, inner sep =2 pt, outer sep= 0 pt, scale=0.2](GH) at (1.5,1) {};
      \node[shape = rectangle,draw, fill=lightgray, text= black, inner sep =2 pt, outer sep= 0 pt, scale=1.7](H) at (2.1,1) {\scriptsize 4};
       \node[shape = rectangle,draw, fill= lightgray, text= black, inner sep =2 pt, outer sep= 0 pt, scale=0.01](I) at (3,1) {};
        \node[shape = rectangle,draw, fill=white, text= black, inner sep =2 pt, outer sep= 0 pt, scale=1.7](DD) at (1.5,0) {\scriptsize 3};
        \node[shape = rectangle,draw, fill=white, text= black, inner sep =2 pt, outer sep= 0 pt, scale=1.7](J) at (0.9,-1) {\scriptsize 2};
     \node[shape = rectangle,draw, fill= black, text= black, inner sep =2 pt, outer sep= 0 pt, scale=0.2](JK) at (1.5,-1) {};
    \node[shape = rectangle,draw, fill= lightgray, text= black, inner sep =2 pt, outer sep= 0 pt, scale=1.7](K) at (2.1,-1) {\scriptsize 5};
     \node[shape = rectangle,draw, fill= lightgray, text= black, inner sep =2 pt, outer sep= 0 pt, scale=0.01](II) at (3,-1) {};
      \node[shape = rectangle,draw, fill= lightgray, text= black, inner sep =2 pt, outer sep= 0 pt, scale=0.01](L) at (3,0) {};
       \node[shape = rectangle,draw, fill= black, text= white, inner sep =2 pt, outer sep= 0 pt, scale=0.2](M) at (3.5,0) {};
        \node[shape = rectangle,draw, fill= black, text= white, inner sep =2 pt, outer sep= 0 pt, scale=0.2](M1) at (3.5,1) {};
         \node[shape = rectangle,draw, fill= black, text=white, inner sep =2 pt, outer sep= 0 pt, scale=0.2](M2) at (3.5,-1) {};
          \node[shape = rectangle,draw, fill= black, text= white, inner sep =2 pt, outer sep= 0 pt, scale=1.7](M22) at (4,-1) {\scriptsize 6};
        \node[shape = rectangle,draw, fill= black, text= white, inner sep =2 pt, outer sep= 0 pt, scale=1.7](N) at (5,-1) {\scriptsize 8};
          \node[shape = rectangle,draw, fill= black, text= white, inner sep =2 pt, outer sep= 0 pt, scale=0.2](N1) at (5.5,-1) {};
           \node[shape = rectangle,draw, fill= black, text= white, inner sep =2 pt, outer sep= 0 pt, scale=0.2](N2) at (5.5,1) {};
            \node[shape = rectangle,draw, fill= black, text= white, inner sep =2 pt, outer sep= 0 pt, scale=0.2](N3) at (5.5,0) {};
             \node[shape = rectangle,draw, fill= black, text= white, inner sep =2 pt, outer sep= 0 pt, scale=0.2](N4) at (6,0) {};
               \node[shape = rectangle,draw, fill= lightgray, text= black, inner sep =2 pt, outer sep= 0 pt, scale=1.7](O1) at (4.5,1) {\scriptsize 7};
     \draw[semithick](B) -- (D);
     \draw[semithick](D) -- (E);
     \draw[semithick](D) -- (F);
      \draw[semithick](E) -- (G);
      \draw[semithick](GH) -- (DD);
     \draw[semithick](JK) -- (DD);
     \draw[semithick](G) -- (H);
     \draw[semithick](H) -- (I);
     \draw[semithick](I) -- (L);
     \draw[semithick](F) -- (J);
     \draw[semithick](J) -- (K);
     \draw[semithick](II) -- (L);
      \draw[semithick](K) -- (II);
     \draw[semithick](L) -- (M);
         \draw[semithick](M) -- (M1);
         \draw[semithick](M) -- (M2);
        \draw[semithick](M2) -- (M22);
       \draw[semithick](M1) -- (O1);
        \draw[semithick](O1) -- (N2);
       \draw[semithick](N2) -- (N3);
         \draw[semithick](M22) -- (N);
         \draw[semithick](N) -- (N1);
         \draw[semithick](N1) -- (N3);
          \draw[semithick](N3) -- (N4);
       \\
  \end{tikzpicture}}
   \caption{\small  System in Example \ref{eg2}}
   \hfill
   \label{fig2}
\end{figure}
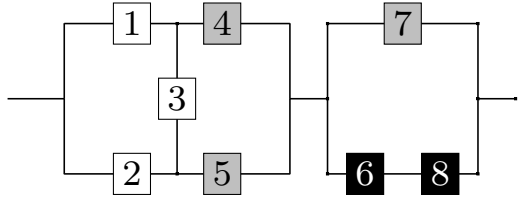

Suppose here that all the components  are independent, where  the components of type $i$ have  common  Weibull reliability functions $\bar{F}_i(t)=e^{-\beta_i t^{\alpha_i}}, \alpha_i, \beta_i>0$ for $i=1,2,3$. In Table \ref{table22} we presented the mean cost rate $Cost_1(v_1, v_2,v_3)$ for given values of $\beta_1=3$, $\beta_2=4,$ $ \beta_3=2$, $\alpha_1=2,$  $\alpha_2=3,$ $ \alpha_3=1$  when the cost parameters are $c_1=1.5, $ $ c_2=1,$ $c_3=2,  c_1^*=0.75$, $c_2^*=0.4,$  $c_3^*=1$ and  $c^{**}=10$. Assume that we have $M_1=7,$  $M_2=4,$ and $M_3=5$ components from types 1, 2, and 3, respectively as spares. Hence, we can choose $v_1\in \{0, 1, 2 \}$, $v_2\in \{0, 1\}$, and $v_3\in \{0,  1, 2\}$ as the redundancy for each type, respectively.

In the left  panel of Table \ref{table22}, the values of $Cost_1$ are computed for different combinations of $v_i$'s. As it is seen,  by adding $v_1=0$, $v_2=1,$ and $ v_3=0$ as the redundant components to groups 1, 2, and 3, respectively, we get the minimum value for the mean cost rate $Cost_1(v_1, v_2,v_3)$.

 Under the assumption that  $\tau$ is the variable of interest, in the right panel of the table,  we have minimized $Cost_2(v_1, v_2,v_3)$ in terms of $\tau$ and have reported  the optimum value of $\tau$,  in the case that the values of $v_i$'s  are kept fixed and known.  It is observed that among all minimized values of $Cost_2$, the least value is obtained for the case that the number of redundant components are $v_1=0$, $v_2=1,$ and $ v_3=0$ for which we have  $\tau=0.375$.
\begin{table}[h!]
\small
\caption{The values of $Cost_1(\mathbf{v})$  and $Cost_2(\mathbf{v})$   in Example \ref{eg2}}\label{table22}
\begin{center}
\begin{tabular}{|ccc|c||c|c|}
\hline
$v_1$ & $v_2$ & $v_3$ & $Cost_1(v_1, v_2,v_3)$ & $\tau_{opt}$ & $Cost_2(v_1, v_2,v_3)$ \\
\hline
0 & 0 & 0 & 39.0424  & { 0.300} & { 29.5929} \\
0 & 1 & 0 & {\bf 37.4142} & {\bf 0.375} & {\bf 28.9959} \\
1 & 0 & 0 & 42.1779 & 0.365 & 34.7858\\
1 & 1 & 0 & 39.0422 & 0.450 & 31.1106 \\
2 & 0 & 0 & 47.4998 & 0.405 & 42.0378 \\
2 & 1 & 0 & 43.0531 & 0.490 & 36.3495 \\
0 & 0 & 1 & 42.2553 & 0.326 & 38.2377 \\
0 & 1 & 1 & 41.2034 & 0.367 & 37.6403 \\
1 & 0 & 1 & 44.1149 & 0.391 & 41.3141 \\
1 & 1 & 1 & 41.9973 & 0.463 & 37.8883 \\
2 & 0 & 1 & 48.5362 & 0.445 & 47.7900 \\
2 & 1 & 1 & 45.5179 & 0.503 & 42.4445 \\
0 & 0 & 2 & 45.5246 & 0.350 & 46.7613 \\
0 & 1 & 2 & 44.8865 & 0.410 & 46.1258 \\
1 & 0 & 2 & 46.1353  & 0.412 & 47.6364 \\
1 & 1 & 2 & 44.8022 & 0.480 & 44.4306 \\
2 & 0 & 2 & 49.7090 & 0.455 & 53.1287 \\
2 & 1 & 2 &  47.7961 & 0.520 & 48.2148 \\
\hline
\end{tabular}
\end{center}
\end{table}
}
\end{Example}
\section{Optimal number of components in series-parallel system}\label{series-parallel}
\begin{figure}[h!]
\small
\centerline{\begin{tikzpicture}[node distance = 3 cm]
  \tikzset{LabelStyle/.style =  {scale=1, fill= white, text=black}}
   \node[shape = circle,draw, fill= white, text= black, inner sep =2 pt, outer sep= 0 pt, scale=0.2](A) at (-3.5,0){} ;
   \node[shape = circle,draw, fill= white, text= black, inner sep =2 pt, outer sep= 0 pt, scale=0.01](D) at (-3,0) {};
   \node[shape = circle,draw, fill= white, text= black, inner sep =2 pt, outer sep= 0 pt, scale=0.01](E) at (-3,1) {};
   \node[shape = circle,draw, fill= white, text= black, inner sep =2 pt, outer sep= 0 pt, scale=0.01](E2) at (-3,2) {};
   \node[shape = circle,draw, fill= white, text= black, inner sep =2 pt, outer sep= 0 pt, scale=0.01](F2) at (-3,-2) {};
   \node[shape = circle,draw, fill= white, text= black, inner sep =2 pt, outer sep= 0 pt, scale=1.7](G) at (-2,1) {\tiny 2};
   \node[shape = circle,draw, fill= white, text= black, inner sep =2 pt, outer sep= 0 pt, scale=1.7](G2) at (-2,2) {\tiny 1};
   \node[shape = circle,draw, fill= white, text= black, inner sep =2 pt, outer sep= 0 pt, scale=0.01](I) at (-1,1) {};
   \node[shape = circle,draw, fill= white, text= black, inner sep =2 pt, outer sep= 0 pt, scale=0.01](I2) at (-1,2) {};
   \node[shape = circle,draw, fill= white, text= black, inner sep =2 pt, outer sep= 0 pt, scale=0.01](L) at (-1,0) {};
   \node[shape = circle,draw, fill=white , text= black, inner sep =2 pt, outer sep= 0 pt, scale=1.7](J) at (-2,-2) {\tiny $n_1$};
   \node[shape = circle,draw, fill= white, text= black, inner sep =2 pt, outer sep= 0 pt, scale=0.01](J2) at (-1,-2) {};
   \node[shape = circle,draw, fill= black, text= black, inner sep =2 pt, outer sep= 0 pt, scale=0.2](M1) at (-2,0) {};
   \node[shape = circle,draw, fill= black, text= black, inner sep =2 pt, outer sep= 0 pt, scale=0.2](M11) at (-2,-0.5) {};
   \node[shape = circle,draw, fill= black, text= black, inner sep =2 pt, outer sep= 0 pt, scale=0.2](M12) at (-2,-1) {};       
     \node[shape = circle,draw, fill= lightgray, text= black, inner sep =2 pt, outer sep= 0 pt, scale=0.01](DD) at (0,0) {};
     \node[shape = circle,draw, fill= lightgray, text= black, inner sep =2 pt, outer sep= 0 pt, scale=0.01](EE) at (0,1) {};
     \node[shape = circle,draw, fill= lightgray, text= black, inner sep =2 pt, outer sep= 0 pt, scale=0.01](EE2) at (0,2) {};
     \node[shape = circle,draw, fill= lightgray, text= black, inner sep =2 pt, outer sep= 0 pt, scale=0.01](FF2) at (0,-2) {};
     \node[shape = circle,draw, fill= lightgray, text= black, inner sep =2 pt, outer sep= 0 pt, scale=1.7](GG) at (1,1) {\tiny 2};
     \node[shape = circle,draw, fill= lightgray, text= black, inner sep =2 pt, outer sep= 0 pt, scale=1.7](GG2) at (1,2) {\tiny 1};
     \node[shape = circle,draw, fill= lightgray, text= black, inner sep =2 pt, outer sep= 0 pt, scale=1.7](JJ) at (1,-2) {\tiny $n_2$};
     \node[shape = circle,draw, fill= lightgray, text= black, inner sep =2 pt, outer sep= 0 pt, scale=0.01](JJ2) at (2,-2) {};
    \node[shape = circle,draw, fill= lightgray, text= black, inner sep =2 pt, outer sep= 0 pt, scale=0.01](DDD) at (2,0) {};
     \node[shape = circle,draw, fill= lightgray, text= black, inner sep =2 pt, outer sep= 0 pt, scale=0.01](II) at (2,1) {};
     \node[shape = circle,draw, fill= lightgray, text= black, inner sep =2 pt, outer sep= 0 pt, scale=0.01](II2) at (2,2) {};
     \node[shape = circle,draw, fill= lightgray, text= black, inner sep =2 pt, outer sep= 0 pt, scale=0.2](M) at (1,0) {};
     \node[shape = circle,draw, fill= lightgray, text= black, inner sep =2 pt, outer sep= 0 pt, scale=0.2](M1) at (1,-0.5) {};
     \node[shape = circle,draw, fill= lightgray, text= black, inner sep =2 pt, outer sep= 0 pt, scale=0.2](M2) at (1,-1) {};
     \node[shape = circle,draw, fill= black, text= white, inner sep =2 pt, outer sep= 0 pt, scale=0.01](D3) at (3,0) {};
     \node[shape = circle,draw, fill= black, text= white, inner sep =2 pt, outer sep= 0 pt, scale=0.01](K) at (3,1) {};
     \node[shape = circle,draw, fill= black, text= white, inner sep =2 pt, outer sep= 0 pt, scale=0.01](B) at (3,2) {};
     \node[shape = circle,draw, fill= black, text= white, inner sep =2 pt, outer sep= 0 pt, scale=0.01](P) at (3,-2) {};
     \node[shape = circle,draw, fill= black, text= white, inner sep =2 pt, outer sep= 0 pt, scale=1.70](O) at (4,1) {\tiny 2};
     \node[shape = circle,draw, fill= black, text= white, inner sep =2 pt, outer sep= 0 pt, scale=1.70](Q) at (4,2) {\tiny 1};
     \node[shape = circle,draw, fill= black, text= white, inner sep =2 pt, outer sep= 0 pt, scale=1.70](R) at (4,-2) {\tiny $n_3$};
     \node[shape = circle,draw, fill= black, text= white, inner sep =2 pt, outer sep= 0 pt, scale=0.01](S) at (5,-2) {};
     \node[shape = circle,draw, fill= black, text= white, inner sep =2 pt, outer sep= 0 pt, scale=0.01](D4) at (5,0) {};
     \node[shape = circle,draw, fill= black, text= white, inner sep =2 pt, outer sep= 0 pt, scale=0.01](T) at (5,1) {};
     \node[shape = circle,draw, fill= black, text= white, inner sep =2 pt, outer sep= 0 pt, scale=0.01](U) at (5,2) {};
     \node[shape = circle,draw, fill= black, text= white, inner sep =2 pt, outer sep= 0 pt, scale=0.20](MM) at (4,0) {};
     \node[shape = circle,draw, fill= black, text= white, inner sep =2 pt, outer sep= 0 pt, scale=0.20](V) at (4,-0.5) {};
     \node[shape = circle,draw, fill= black, text= white, inner sep =2 pt, outer sep= 0 pt, scale=0.20](W) at (4,-1) {};
     \node[shape = circle,draw, fill= black, text= white, inner sep =2 pt, outer sep= 0 pt, scale=0.20](M3) at (5.5,0) {};
     \draw[semithick](A) -- (D);
     \draw[semithick](D) -- (E);
     \draw[semithick](E) -- (E2);
     \draw[semithick](E) -- (G);
     \draw[semithick](G) -- (I);
     \draw[semithick](G2) -- (I2);
     \draw[semithick](E2) -- (G2);
     \draw[semithick](D) -- (F2);
     \draw[semithick](F2) -- (J);
     \draw[semithick](J) -- (J2);
     \draw[semithick](J2) -- (L);
     \draw[semithick](I2) -- (L);
     \draw[semithick](I) -- (I2);
     \draw[semithick](I) -- (L);
     \draw[semithick](L) -- (DD);
     \draw[semithick](DD) -- (EE);
     \draw[semithick](EE) -- (EE2);
     \draw[semithick](DD) -- (FF2);
     \draw[semithick](EE) -- (GG);
     \draw[semithick](EE2) -- (GG2);
     \draw[semithick](FF2) -- (JJ);
     \draw[semithick](JJ) -- (JJ2);
     \draw[semithick](GG) -- (II);
     \draw[semithick](GG2) -- (II2);
     \draw[semithick](II) -- (II2);
     \draw[semithick](II2) -- (DDD);
     \draw[semithick](JJ2) -- (DDD);
     \draw[semithick](D3) -- (K);
     \draw[semithick](K) -- (B);
     \draw[semithick](D3) -- (P);
     \draw[semithick](K) -- (O);
     \draw[semithick](B) -- (Q);
     \draw[semithick](P) -- (R);
     \draw[semithick](R) -- (S);
     \draw[semithick](O) -- (T);
     \draw[semithick](Q) -- (U);
     \draw[semithick](T) -- (U);
     \draw[semithick](T) -- (D4);
     \draw[semithick](S) -- (D4);
     \draw[semithick](DDD) -- (D3);
     \draw[semithick](D4) -- (M3);
          \end{tikzpicture}}
   \caption{\small  A series-parallel system with 3 subsystems }
   \hfill
   \label{figsp}
\end{figure}
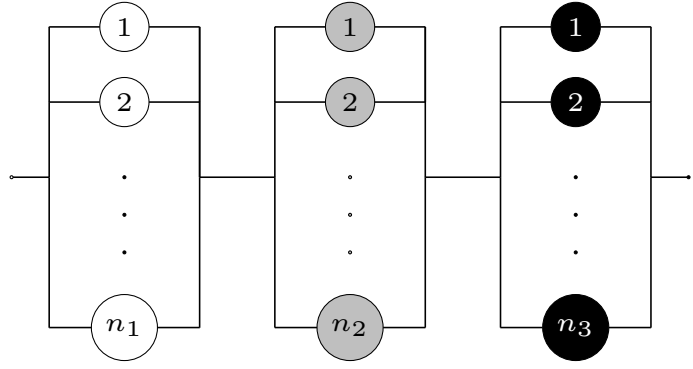

\quad An important subclass of coherent systems is the class of series-parallel systems. A series-parallel system is a series structure of $L$ parallel subsystems, $L\geq 1$; see e.g. Figure \ref{figsp}.
The purpose here is to find the optimal number of the components in the $l$th parallel subsystem, under the condition that there are available $M_l$ components of type $l$, where the components in $l$th subsystem are exchangeable dependent having common reliability function $\bar{F}_l$, $l=1,\ldots, L$. Furthermore, suppose that the random failure times of the components of different types are dependent. The dependency structure in the system is built with a copula function $\hat{C}$, as described in Section 2.
Under the mean cost rate criteria defined in Section 2, the problem of optimal allocation is to find the optimal values of $n_l$ for each subsystem so that (\ref{cost1}) or (\ref{cost2}) is minimized.

In the following, we  provide the corresponding expressions for the cost-based functions in a  series-parallel system.
First note that for this system, we have
\begin{align*}
\Phi(l_1,\ldots, l_L)=\left\{
\begin{tabular}{ll}
1 & $\forall j\in\{1,\ldots,L\}: l_j\geq 1$\\
0 & o.w.\\
\end{tabular}
\right.
\end{align*}
Hence from \eqref{Fbar-dep}, we get
\begin{align*}
\bar{F}_T(t)&={\sum_{l_1=1}^{n_1}\ldots\sum_{l_L=1}^{n_L}}{\sum_{i_1=0}^{n_1-l_1}\ldots\sum_{i_L=0}^{n_L-l_L}}(-1)^{i_1+\ldots+ i_L}\binom{n_1}{l_1}\ldots \binom{n_L}{l_L}\binom{n_1-l_1}{i_1}\ldots \binom{n_L-l_L}{i_L}\nonumber\\
&\times \hat{C}(\underbrace{\bar{F}_1(t)}_{i_1+l_1},\underbrace{1}_{n_1-(i_1+l_1)},\ldots,\underbrace{\bar{F}_L(t)}_{i_L+l_L},\underbrace{1}_{n_L-(i_L+l_L). }),
\end{align*}
and in the especial case of independent components, we derive  from \eqref{ind}
\begin{align*}
\bar{F}_T(t)&=\prod_{l=1}^L (1-[1-\bar{F}_l(t)]^{n_l}).
\end{align*}
\subsection*{Cost function at the  system failure}
\quad In a similar manner to Subsection \ref{cost}, the mean cost rate function for system failure is defined as:
\begin{align}\label{cost78}
Cost_3(\mathbf{n})=\frac{\sum_{i=1}^Lc_iE(X_i(T))+\sum_{i=1}^L c_i^*E(n_i-X_i(T))+c^{**}}{E(T)}
\end{align}
where $\mathbf{n}=(n_1,\ldots, n_L)$, and
\begin{align}\label{3}
E(X_i(T))&=n_i \sum_{m_1=1}^{n_1}\ldots\sum_{m_i=0}^{n_i-1}\ldots\sum_{m_L=1}^{n_L}\binom{n_1}{m_1}\ldots \binom{n_i-1}{m_i}\ldots\binom{n_L}{m_L}\int_{0}^{\infty}\lim_{\delta\rightarrow 0}\frac{A_{\mathbf{m}}^{(i)}(t,\delta)}{\delta}dt,
\end{align}
 in which $A_{\mathbf{m}}^{(i)}(t,\delta)$ is introduced in \eqref{Adelta}. For the independent components \eqref{3} reduces to the following expression
\begin{small}
\begin{align}\label{78}
E(X_i(T))
&=n_i \int_0^{\infty} \prod_{l=1,l\neq i}^L (1-(1-\bar{F}_l(t))^{n_l}) dF_{i}(t).
\end{align}
\end{small}
If $L=1$ then the system becomes a parallel system with $n_1$ components, and hence in this case $E(X_1(T))=n_1$.

{
\quad For the considered series-parallel system in which the components of subsystems are independent,  \cite{Eryilmazetal2020} gained a similar result for $E(X_i(T))$  in \eqref{78}. Subsequently, they found the optimal numbers of components in each subsystem based on the minimization of cost function \eqref{cost78} under the constraints on the total allotted cost for replacing failed components and the total allotted cost for rejuvenation of unfailed ones. Hence, our results in this subsection   may be considered as an extension of  their work to the case of dependent components.} { Also, \cite{Dembinska2021} discussed the similar problem for the case that  the lifetime distributions of components are discrete;   especially they obtained some results for discrete phase-type distribution.}
\subsection*{Cost function based on preventive  replacement}
The mean cost rate function of the system for age replacement at time $\min(\tau, T)$ is defined as
\begin{align*}
Cost_4(\mathbf{n})=\frac{M_1(\mathbf{n})P(T\leq \tau)+M_2(\mathbf{n})P(T> \tau)}{E(\min(\tau,T))}
\end{align*}
where
\begin{align*}
M_1(\mathbf{n})=\sum_{i=1}^L c_i E(X_i(T)|T\leq \tau)+\sum_{i=1}^L c^*_i E(n_i-X_i{(T)}|T\leq \tau)+c^{**},
\end{align*}
and
\begin{align*}
M_2(\mathbf{n})=\sum_{i=1}^L c_i E(N_i{(\tau)}|T>\tau)+\sum_{i=1}^L c^*_i E(n_i-N_i{(\tau)}|T>\tau).
\end{align*}
Using the formula for the survival signature of the series-parallel system, from the results given in Subsection 2.2, we get the following expressions.
\begin{align*}
E(N_i(\tau)|T>\tau)=\frac{1}{\bar{F}(\tau)}\sum_{j_1=0}^{n_1-1}\ldots\sum_{j_L=0}^{n_L-1} j_i \binom{n_1}{j_1}\ldots \binom{n_L}{j_L}B(\tau, j_1,\dots,j_L)
\end{align*}
and { for $n_i\geq 2$}
\begin{align}\label{iop}
&E(X_i(T)|T\leq \tau)=\frac{n_i}{1-\bar{F}_T(\tau)} \Bigg[\sum_{m_1=1}^{n_1}\ldots\sum_{m_i=1}^{n_i-1}\ldots\sum_{m_L=1}^{n_L}
{\sum_{l_1=0}^{m_1}\ldots\sum_{l_L=0}^{m_L}}\bigg[\prod_{j=1,j\neq i}^L\binom{n_j}{m_j}\binom{m_j}{l_j}\bigg]\nonumber\\
&\times\binom{n_i-1}{m_i}\binom{m_i}{l_i}\int_{0}^{\tau}\lim_{\delta\rightarrow 0}\frac{1}{\delta} A_{\mathbf{m}, \mathbf{l}}^{(i)}(s,s+\delta, \tau)ds -
 \sum_{m_1=1}^{n_1}\ldots\sum_{m_i=1}^{n_i-1}\ldots\sum_{m_L=1}^{n_L}
{\sum_{l_1=1}^{m_1}\ldots\sum_{l_L=1}^{m_L}}\nonumber\\
&\bigg[\prod_{j=1,j\neq i}^L\binom{n_j}{m_j}\binom{m_j}{l_j}\bigg]\binom{n_i-1}{m_i}\binom{m_i}{l_i}\int_{0}^{\tau}\lim_{\delta\rightarrow 0}\frac{1}{\delta} A_{\mathbf{m}, \mathbf{l}}^{(i)}(s,s+\delta, \tau)ds \Bigg].
\end{align}
If $n_i=1$ then it is easily deduced that $E(X_i(T)|T\leq \tau)=\frac{{F}_i(\tau)}{1-\bar{F}_T(\tau)}$.

\vspace{0.5cm}

\begin{Example}\label{ex-series}
{\rm Consider a series-parallel system with $L=3$ subsystems  and assume that there are $M_1=2, M_2=3$ and $M_3=3$ components from types $1, 2$, and $3$, respectively, to construct the system. Suppose that the joint reliability function  of the components lifetimes follow the multivariate Pareto model given by
 \begin{align*}
&P\big( T_1^{(1)}>t_1^{(1)}, \ldots,  T_{n_1}^{(1)}>t_{n_1}^{(1)}, \ldots,  T_1^{(L)}>t_L^{(L)}, \ldots,  T_{n_L}^{(L)}>t_{n_L}^{(L)}\big)\\
 &\hspace{2cm}=\left[ 1+\theta_1\sum_{i=1}^{n_1}t_i^{(1)}+\ldots+\theta_L\sum_{i=1}^{n_L}t_i^{(L)}     \right]^{-\alpha} ,
 \end{align*}
for $\theta_i>0, i=1, \ldots, L$, and $\alpha>0$.
In fact, the corresponding survival copula is
\begin{align*}
\hat{C}(u_1,\ldots, u_n)=\left(u_1^{-1/\alpha}+\ldots + u_n^{-1/\alpha}-n+1  \right)^{-\alpha},
\end{align*}
and the marginal reliability functions of the components in the subsystems are $\bar{F}_i(t)=(1+\theta_i t)^{-\alpha}$, $i=1,2, \ldots, L$.
{

First note that for the described system, we have
\begin{small}
\begin{align*}
\int_{0}^{\infty}\lim_{\delta\rightarrow 0}\frac{A_{\mathbf{m}}^{(i)}(t,\delta)}{\delta}dt&={\sum_{j_1=0}^{n_1-m_1}\ldots\sum_{j_i=0}^{n_i-m_i-1}\ldots\sum_{j_L=0}^{n_L-m_L}}(-1)^{j_1+\ldots + j_L}\binom{n_1-m_1}{j_1}\ldots \binom{n_i-m_i-1}{j_i}\ldots\binom{n_L-m_L}{j_L}\\
&\times \frac{\theta_i}{\theta_i(m_i+j_i+1)+\sum_{l=1, l\neq i}^{L}\theta_l(m_l+j_l)}
\end{align*}
\end{small}
By replacing these expressions  in \eqref{3}, $E(X_i(T)), i=1, 2, 3$  are obtained. 
}

Let ${\mathbf{\theta}}=(0.4, 0.2, 0.3)$,  $\mathbf{c}=(1.5, 2, 3)$, $\mathbf{c}^*=(0.3, 0.75, 1)$, $c^{**}=8$, and $\alpha=2$. The values of the mean cost rate function $Cost_3$ are obtained for all combinations of $n_1, n_2$, and $n_3$, such that $n_1\in\{1, 2\}, n_2\in\{1, 2, 3\}$, and $n_3\in\{1,2,3\}$. Also, under the age replacement strategy in $\tau=1$, the values of mean cost rate function $Cost_4$ are calculated for different values  $n_1, n_2$ and $n_3$. 

{
For computing $Cost_4(\tau)$, we use the following simplified expressions:
\begin{align*}
&E(N_i(\tau)|T>\tau)=\frac{1}{\bar{F}_T(\tau)}\sum_{j_1=0}^{n_1-1}\sum_{j_2=0}^{n_2-1}\sum_{j_3=0}^{n_3-1} j_i \binom{n_1}{j_1} \binom{n_2}{j_2} \binom{n_3}{j_3}{\sum_{b_1=0}^{j_1}\sum_{b_2=0}^{j_2}\sum_{b_3=0}^{j_3}}(-1)^{b_1+ b_2+ b_L}\\
&\binom{j_1}{b_1}\binom{j_2}{b_2}\binom{j_3}{b_3}[1+\theta_1\tau(n_1-j_1+b_1)+\theta_2\tau(n_2-j_2+b_2)+\theta_3\tau(n_3-j_3+b_3)]^{-\alpha},~i=1, 2, 3,
\end{align*}
where
\begin{align*}
\bar{F}_T(\tau)&={\sum_{l_1=1}^{n_1}\sum_{l_2=1}^{n_2}\sum_{l_3=1}^{n_3}}{\sum_{i_1=0}^{n_1-l_1}\sum_{i_2=0}^{n_2-l_2}\sum_{i_3=0}^{n_3-l_3}}(-1)^{i_1+ i_2+ i_3}\binom{n_1}{l_1}\binom{n_2}{l_2}\binom{n_3}{l_3}\binom{n_1-l_1}{i_1}\binom{n_2-l_2}{i_2} \binom{n_3-l_3}{i_3}\nonumber\\
&\big(1+\theta_1 \tau(i_1+l_1)+\theta_2 \tau(i_2+l_2)+\theta_3 \tau(i_3+l_3)\big)^{-\alpha}.
\end{align*}


Also, we obtain  $E(X_i(T)|T\leq \tau),  i=1, 2, 3$ by placing the following quantity in \eqref{iop},
\begin{small}
\begin{align*}
&\int_{0}^{\tau}\lim_{\delta\rightarrow 0}\frac{1}{\delta} A_{\mathbf{m}, \mathbf{l}}^{(i)}(s,s+\delta, \tau)ds={\sum_{j_1=0}^{n_1-m_1}\ldots\sum_{j_i=0}^{n_i-m_i-1}\ldots\sum_{j_L=0}^{n_L-m_L}}(-1)^{j_1+\ldots + j_L}\binom{n_1-m_1}{j_1}\ldots \binom{n_i-m_i-1}{j_i}\ldots\binom{n_L-m_L}{j_L}\\
&\times {\sum_{d_1=0}^{m_1-l_1}\ldots\sum_{d_i=0}^{m_i-l_i}\ldots\sum_{d_L=0}^{m_L-l_L}}(-1)^{d_1+\ldots + d_L}\binom{m_1-l_1}{d_1}\ldots \binom{m_L-l_L}{d_L}\\
&\times \frac{\theta_i[\big(1+\tau(\sum_{k=1}^L \theta_k(l_k+d_k)) \big)^{-\alpha}-\big(1+\tau(\theta_i(m_i+j_i+1)+\sum_{k=1, \neq i}^L \theta_k(m_k+j_k))\big)^{-\alpha}]}{\theta_i(m_i-l_i+j_i-d_i+1)+\sum_{k=1, \neq i}^L \theta_k(m_k-l_k+j_k-d_k)}.
\end{align*}
\end{small}
}

The results are given in Table \ref{table222}.
It is seen from the results that based on objective function $Cost_3(.)$ the optimal series-parallel system has $n_1=2, n_2=3, n_3=2$ components. Since the reliability of the second type is greater than the other two types, it is expected that more components for the second subsystem lead to a reduction in the mean cost rate of system failure. 
Also, $ n_1 = 2, n_2 = 2, n_3 = 2 $ are the optimal number of components in the subsystems so as to minimize the average cost rate of the age replacement policy.
\begin{table}[t]
\small
\caption{The values of $Cost_3(\mathbf{n})$  and $Cost_4(\mathbf{n})$   in Example \ref{ex-series}}\label{table222}
\begin{center}
\begin{tabular}{|ccc|c|c|}
\hline

$n_1$ & $n_2$ & $n_3$ & {$Cost_3(n_1,n_2,n_3)$} & $Cost_4(n_1,n_2,n_3)$ \\

\hline
1 & 1 & 1 & 10.3750 & 18.2784  \\
1 & 1 & 2 & 9.6460 & 18.6515 \\
1 & 1 & 3 & 10.0967 & 21.2331 \\

1 & 2 & 1 & 9.7277 & 18.1732 \\
1 & 2 & 2 & 8.8857 & 18.1463  \\
1 & 2 & 3 & 9.1886 & 20.4388 \\

1 & 3 & 1 & 10.0842 & 19.8532  \\
1 & 3 & 2 & 9.0989 & 19.5615 \\
1 & 3 & 3 & 9.3157 & 21.7702\\

2 & 1 & 1 & 8.7073  & 16.2100\\
2 & 1 & 2 & 7.9885 & 16.1742\\
2 & 1 & 3 & 8.2879 & 18.3280  \\

2 & 2 & 1 & 8.0593  & 15.8127\\
2 & 2 & 2 & 7.4747 & {\bf 12.2278}\\
2 & 2 & 3 & 7.4835 & {13.5080} \\

2 & 3 & 1 & 8.2888 & 17.1959 \\
2 & 3 & 2 & {\bf 7.4068} & 13.0874\\
2 & 3 & 3 & 7.5279 & {14.2977 }\\
\hline

\end{tabular}
\end{center}
\end{table}

}
\end{Example}
\section{Conclusions}
\quad In this paper, we studied the optimal number of redundancy allocation in an $n$-component coherent system consists of heterogeneous components. We assumed that the system has been built up of $L$ different types of components, $L\geq 1$, where there are $n_i$ components of type $i$ and $\sum_{i=1}^{L}n_i=n$. We assumed that the components of the different types in the system are statistically dependent. The system reliability function was modeled by the notion of survival signature in terms of a given survival copula function. We further assumed $M_i$ components available as spares for the components  of type $i$. We investigated the number of active redundant components $v_i$, $n_iv_i\leq M_i$,  that can be added to each component of type $i$ such that  the imposed cost functions are minimized, $i=1,\dots, L$.  { We first proposed a cost function in terms of the costs of renewing the failed components and the costs of refreshing the alive components at the time of the system failure. In the sequel, we proposed an  cost-based function in terms of the costs of the renewing (refreshing) the failed (alive) components at the system failure time or at a predetermined time $\tau$, whichever occurs first.}  In the last part of the paper, we studied under the settings of the first part, the particular case that the system is a series-parallel system. We derived the formulas for the proposed cost functions and using them investigated the optimal number of the components in each parallel subsystem. The expressions for the proposed cost functions were derived using the mixture representation of the system reliability function based on the notion of survival signature.  The results were examined numerically for some particular coherent systems.
The proposed mean cost rate functions simultaneously consider the cost of the system and its $MTTF$  (which is directly related to its reliability). Hence this optimization problem can be viewed as a bi-objective reliability-redundancy allocation problem but with  a more comfortable setup.


  In this study, we considered the general case that the components of the same group are exchangeable and the components of different groups and dependent.  Although these assumptions are more realistic and hence increase the range of applications of our results,  however, obviously lead to the complexity of the formulas. Even a more realistic case is the situation that the components in each group are dependent in a more general sense than that of exchangeability. Developing results in this direction may be considered as a future study.
Here, we presumed the active redundancy for components. Allocating the other variants of spares, i.e., cold and warm standby, for coherent systems may be investigated as some interesting problems for future studies.

%
%
%

\section*{Acknowledgments:} { We would like to thank Associate Editor and two anonymous reviewers  for their constructive comments and suggestions which led to improvements of the presentation of the paper.} Asadi's research work was performed in IPM Isfahan branch and was in part supported by a grant from IPM (No. 1400620212).

\section*{Appendix}

{\bf Proof of Theorem \ref{Am}}

Define the events $M$, $N^c$  and $L^C$ as follows
\begin{align*}
M&\equiv\{{T_1^{(1)}>t},\ldots,{T_{m_1}^{(1)}>t},{T_2^{(i)}>t},\ldots,{T_{m_i+1}^{(i)}>t},{T_1^{(L)}>t},\ldots,{T_{m_L}^{(L)}>t} \}\\
N^c &\equiv\{{T_{m_1+1}^{(1)}<t},\ldots, {T_{n_1}^{(1)}<t}, {T_{m_i+2}^{(i)}<t},\ldots, {T_{n_i}^{(i)}<t}, {T_{m_L+1}^{(L)}<t},\ldots, {T_{n_L}^{(L)}<t}   \}\\
L^c &\equiv\{t<T_1^{(i)}<t+\delta\}.
\end{align*}
Then $A_{\mathbf{m}}(t, \delta)$ equals to the following
\begin{align}\label{av3}
A_{\mathbf{m}}^{(i)}(t, \delta)&=P(M\cap N^c \cap L^c)=P(M)-P(M\cap N)-P(M \cap L)+P(M \cap N \cap L).
\end{align}
Evidently, we have
\begin{align*}
&P(M)=\hat{C}(\underbrace{\bar{F}_1(t)}_{m_1},\underbrace{1}_{n_1-m_1},\ldots, \underbrace{\bar{F}_i(t)}_{m_i+1},\underbrace{1}_{n_i-m_i-1}, \ldots, \underbrace{\bar{F}_L(t)}_{m_L},\underbrace{1}_{n_L-m_L})
\end{align*}
and
\begin{align*}
&P(M \cap L)=\hat{C}(\underbrace{\bar{F}_1(t)}_{m_1},\underbrace{1}_{n_1-m_1},\ldots, \underbrace{\bar{F}_i(t)}_{m_i}, \bar{F}_i(t+\delta),\underbrace{1}_{n_i-m_i-1}, \ldots, \underbrace{\bar{F}_L(t)}_{m_L},\underbrace{1}_{n_L-m_L}).
\end{align*}
Note that we can write
\begin{align*}
N=\cup_{j=m_1+1}^{n_1} \{T_{j}^{(1)}>t\}\cup \ldots \cup_{j=m_i+2}^{n_i} \{T_{j}^{(i)}>t\}\cup \ldots  \cup_{j=m_L+1}^{n_L} \{T_{j}^{(L)}>t\}.
\end{align*}
Therefore, we can easily see that
\begin{small}
\begin{align*}
&P(M\cap N)\\
&=\sum_{l=1}^{n-\sum_{i=1}^L m_i -1}(-1)^{l+1}\underset{j_1+\ldots +j_L=l}{\sum_{j_1=0}^{n_1-m_1}\ldots \sum_{j_i=0}^{n_i-m_i-1}\ldots \sum_{j_L=0}^{n_L-m_L}}
\binom{n_1-m_1}{j_1}\ldots \binom{n_i-m_i-1}{j_i}\ldots\binom{n_L-m_L}{j_L}\\
&\ \ \ \times  \hat{C}(\underbrace{\bar{F}_1(t)}_{m_1+j_1},\underbrace{1}_{n_1-m_1-j_1},\ldots, \underbrace{\bar{F}_i(t)}_{m_i+j_i+1},\underbrace{1}_{n_i-m_i-j_i-1}, \ldots, \underbrace{\bar{F}_L(t)}_{m_L+j_L},\underbrace{1}_{n_L-m_L-j_L}).
\end{align*}
\end{small}
If we subtract  $P(M)$ from both sides of this equation, then we have

\begin{small}
\begin{align}\label{av2}
&P(M\cap N)-P(M)\nonumber\\
&=\sum_{l=1}^{n-\sum_{i=1}^L m_i -1}(-1)^{l+1}\underset{j_1+\ldots +j_L=l}{\sum_{j_1=0}^{n_1-m_1}\ldots \sum_{j_i=0}^{n_i-m_i-1}\ldots \sum_{j_L=0}^{n_L-m_L}}\binom{n_1-m_1}{j_1}\ldots \binom{n_i-m_i-1}{j_i}\ldots\binom{n_L-m_L}{j_L}\nonumber\\
&\ \ \ \times  \hat{C}(\underbrace{\bar{F}_1(t)}_{m_1+j_1},\underbrace{1}_{n_1-m_1-j_1},\ldots, \underbrace{\bar{F}_i(t)}_{m_i+j_i+1},\underbrace{1}_{n_i-m_i-j_i-1}, \ldots, \underbrace{\bar{F}_L(t)}_{m_L+j_L},\underbrace{1}_{n_L-m_L-j_L})\nonumber\\
&\hspace{7.5cm}-[\hat{C}(\underbrace{\bar{F}_1(t)}_{m_1},\underbrace{1}_{n_1-m_1},\ldots, \underbrace{\bar{F}_i(t)}_{m_i+1},\underbrace{1}_{n_i-m_i-1}, \ldots, \underbrace{\bar{F}_L(t)}_{m_L},\underbrace{1}_{n_L-m_L})]\nonumber\\
&=\sum_{l=0}^{n-\sum_{i=1}^L m_i -1}(-1)^{l+1}\underset{j_1+\ldots +j_L=l}{\sum_{j_1=0}^{n_1-m_1}\ldots \sum_{j_i=0}^{n_i-m_i-1}\ldots \sum_{j_L=0}^{n_L-m_L}}\binom{n_1-m_1}{j_1}\ldots \binom{n_i-m_i-1}{j_i}\ldots\binom{n_L-m_L}{j_L}\nonumber\\
&\ \ \ \times  \hat{C}(\underbrace{\bar{F}_1(t)}_{m_1+j_1},\underbrace{1}_{n_1-m_1-j_1},\ldots, \underbrace{\bar{F}_i(t)}_{m_i+j_i+1},\underbrace{1}_{n_i-m_i-j_i-1}, \ldots, \underbrace{\bar{F}_L(t)}_{m_L+j_L},\underbrace{1}_{n_L-m_L-j_L})\nonumber\\
&={\sum_{j_1=0}^{n_1-m_1}\ldots \sum_{j_i=0}^{n_i-m_i-1}\ldots \sum_{j_L=0}^{n_L-m_L}}(-1)^{j_1+\ldots +j_L+1}\binom{n_1-m_1}{j_1}\ldots \binom{n_i-m_i-1}{j_i}\ldots\binom{n_L-m_L}{j_L}\nonumber\\
&\ \ \ \times  \hat{C}(\underbrace{\bar{F}_1(t)}_{m_1+j_1},\underbrace{1}_{n_1-m_1-j_1},\ldots, \underbrace{\bar{F}_i(t)}_{m_i+j_i+1},\underbrace{1}_{n_i-m_i-j_i-1}, \ldots, \underbrace{\bar{F}_L(t)}_{m_L+j_L},\underbrace{1}_{n_L-m_L-j_L}).
\end{align}
\end{small}
Similarly, we have
\begin{small}
\begin{align*}
&P(M\cap N \cap L)\nonumber\\
&=\sum_{l=1}^{n-\sum_{i=1}^L m_i -1}(-1)^{l+1}\underset{j_1+\ldots +j_L=l}{\sum_{j_1=0}^{n_1-m_1}\ldots \sum_{j_i=0}^{n_i-m_i-1}\ldots \sum_{j_L=0}^{n_L-m_L}}\binom{n_1-m_1}{j_1}\ldots \binom{n_i-m_i-1}{j_i}\ldots\binom{n_L-m_L}{j_L}\nonumber\\
&\ \ \ \times \hat{C}(\underbrace{\bar{F}_1(t)}_{m_1+j_1},\underbrace{1}_{n_1-m_1-j_1},\ldots, \underbrace{\bar{F}_i(t)}_{m_i+j_i}, \bar{F}_i(t+\delta),\underbrace{1}_{n_i-m_i-j_i-1}, \ldots, \underbrace{\bar{F}_L(t)}_{m_L+j_L},\underbrace{1}_{n_L-m_L-j_L}).\nonumber \\
\end{align*}
\begin{align}\label{av1}
&P(M\cap N \cap L)-P(M\cap L)\nonumber\\
&=\sum_{l=0}^{n-\sum_{i=1}^L m_i -1}(-1)^{l+1}\underset{j_1+\ldots +j_L=l}{\sum_{j_1=0}^{n_1-m_1}\ldots \sum_{j_i=0}^{n_i-m_i-1}\ldots \sum_{j_L=0}^{n_L-m_L}}\binom{n_1-m_1}{j_1}\ldots \binom{n_i-m_i-1}{j_i}\ldots\binom{n_L-m_L}{j_L}\nonumber\\
&\ \ \ \times \hat{C}(\underbrace{\bar{F}_1(t)}_{m_1+j_1},\underbrace{1}_{n_1-m_1-j_1},\ldots, \underbrace{\bar{F}_i(t)}_{m_i+j_i}, \bar{F}_i(t+\delta),\underbrace{1}_{n_i-m_i-j_i-1}, \ldots, \underbrace{\bar{F}_L(t)}_{m_L+j_L},\underbrace{1}_{n_L-m_L-j_L})\nonumber\\
&={\sum_{j_1=0}^{n_1-m_1}\ldots \sum_{j_i=0}^{n_i-m_i-1}\ldots \sum_{j_L=0}^{n_L-m_L}}(-1)^{j_1+\ldots +j_L+1}\binom{n_1-m_1}{j_1}\ldots \binom{n_i-m_i-1}{j_i}\ldots\binom{n_L-m_L}{j_L}\nonumber\\
&\ \ \ \times \hat{C}(\underbrace{\bar{F}_1(t)}_{m_1+j_1},\underbrace{1}_{n_1-m_1-j_1},\ldots, \underbrace{\bar{F}_i(t)}_{m_i+j_i}, \bar{F}_i(t+\delta),\underbrace{1}_{n_i-m_i-j_i-1}, \ldots, \underbrace{\bar{F}_L(t)}_{m_L+j_L},\underbrace{1}_{n_L-m_L-j_L}).
\end{align}
\end{small}
Then replacing \eqref{av2} and \eqref{av1} in \eqref{av3} we have
\begin{small}
\begin{align*}
&A_{\mathbf{m}}^{(i)}(t, \delta)=
-{\sum_{j_1=0}^{n_1-m_1}\ldots \sum_{j_i=0}^{n_i-m_i-1}\ldots \sum_{j_L=0}^{n_L-m_L}}
(-1)^{j_1+\ldots +j_L+1}\binom{n_1-m_1}{j_1}\ldots \binom{n_i-m_i-1}{j_i}\ldots\binom{n_L-m_L}{j_L}\\
&\hspace{5cm} \times  \hat{C}(\underbrace{\bar{F}_1(t)}_{m_1+j_1},\underbrace{1}_{n_1-m_1-j_1},\ldots, \underbrace{\bar{F}_i(t)}_{m_i+j_i+1},\underbrace{1}_{n_i-m_i-j_i-1}, \ldots, \underbrace{\bar{F}_L(t)}_{m_L+j_L},\underbrace{1}_{n_L-m_L-j_L})\\
&+{\sum_{j_1=0}^{n_1-m_1}\ldots \sum_{j_i=0}^{n_i-m_i-1}\ldots \sum_{j_L=0}^{n_L-m_L}}
(-1)^{j_1+\ldots +j_L+1}\binom{n_1-m_1}{j_1}\ldots \binom{n_i-m_i-1}{j_i}\ldots\binom{n_L-m_L}{j_L}\\
&\hspace{4cm}\times  \hat{C}(\underbrace{\bar{F}_1(t)}_{m_1+j_1},\underbrace{1}_{n_1-m_1-j_1},\ldots, \underbrace{\bar{F}_i(t)}_{m_i+j_i}, \bar{F}_i(t+\delta),\underbrace{1}_{n_i-m_i-j_i-1}, \ldots, \underbrace{\bar{F}_L(t)}_{m_L+j_L},\underbrace{1}_{n_L-m_L-j_L})\\
&={\sum_{j_1=0}^{n_1-m_1}\ldots \sum_{j_i=0}^{n_i-m_i-1}\ldots \sum_{j_L=0}^{n_L-m_L}}
(-1)^{j_1+\ldots +j_L}\binom{n_1-m_1}{j_1}\ldots \binom{n_i-m_i-1}{j_i}\ldots\binom{n_L-m_L}{j_L}\\
&\times \bigg[ \hat{C}(\underbrace{\bar{F}_1(t)}_{m_1+j_1},\underbrace{1}_{n_1-m_1-j_1},\ldots, \underbrace{\bar{F}_i(t)}_{m_i+j_i+1},\underbrace{1}_{n_i-m_i-j_i-1}, \ldots, \underbrace{\bar{F}_L(t)}_{m_L+j_L},\underbrace{1}_{n_L-m_L-j_L})\\
&\hspace{4cm} -\hat{C}(\underbrace{\bar{F}_1(t)}_{m_1+j_1},\underbrace{1}_{n_1-m_1-j_1},\ldots, \underbrace{\bar{F}_i(t)}_{m_i+j_i}, \bar{F}_i(t+\delta),\underbrace{1}_{n_i-m_i-j_i-1}, \ldots, \underbrace{\bar{F}_L(t)}_{m_L+j_L},\underbrace{1}_{n_L-m_L-j_L})\bigg]
\end{align*}
\end{small}
\qed

{\bf Proof of Theorem \ref{Aml}}:\\

Let us define the following events
\begin{small}
\begin{align*}
M&\equiv\{T_1^{(j)}>\tau, \ldots, T_{l_j}^{(j)}>\tau, T_{l_j+1}^{(j)}>s, \ldots, T_{m_j}^{(j)}>s,  \text{for $j=1,...,L, j\neq i$}\\
&\hspace{4.3cm},T_1^{(i)}>\tau, \ldots, T_{l_i}^{(i)}>\tau, T_{l_i+1}^{(i)}>s, \ldots, T_{m_i+1}^{(i)}>s\}\\
N^c&\equiv\{ T_{m_j+1}^{(j)}<s,\ldots, T_{n_j}^{(j)}<s,  \text{for $j=1,...,L, j\neq i$}, T_{m_i+2}^{(i)}<s,\ldots, T_{n_i}^{(i)}<s  \}\\
L^c&\equiv\{T_{l_j+1}^{(j)}<\tau, \ldots, T_{m_j}^{(j)}<\tau,\text{for $j=1,...,L$} \}\\
K^c&\equiv\{T_{m_i+1}^{(i)}<s+\delta \}.
\end{align*}
\end{small}
Hence,
\begin{small}
\begin{align}\label{av5}
A_{\mathbf{m}, \mathbf{l}}^{(i)}(s,s+\delta, \tau)&=P(M\cap N^c \cap L^c \cap K^c)\nonumber\\
&=P(M)-P(M\cap N)-P(M\cap L)-P(M \cap K)+P(M \cap N \cap L)+P(M \cap N \cap K)\nonumber\\
&+P(M \cap K \cap L)-P(M \cap N \cap L \cap K).
\end{align}
\end{small}
It can be easily shown that
\begin{small}
\begin{align*}
P(M)&= \hat{C}(\underbrace{\bar{F}_j(\tau)}_{l_j},\underbrace{\bar{F}_j(s)}_{m_j-l_j},\underbrace{1}_{n_j-m_j},  \text{for $j=1,...,L, j\neq i$}, \underbrace{\bar{F}_i(\tau)}_{l_i}, \underbrace{\bar{F}_i(s)}_{m_i-l_i+1},
\underbrace{1}_{n_i-m_i-1}),\\
P(M \cap K)&=P(M\cap \left\{ T_{m_i+1}^{(i)}>s+\delta \}  \right\})\\
&=\hat{C}(\underbrace{\bar{F}_j(\tau)}_{l_j},\underbrace{\bar{F}_j(s)}_{m_j-l_j},\underbrace{1}_{n_j-m_j},  \text{for $j=1,...,L, j\neq i$}, \underbrace{\bar{F}_i(\tau)}_{l_i}, \underbrace{\bar{F}_i(s)}_{m_i-l_i}, \bar{F}_i(s+\delta),\underbrace{1}_{n_i-m_i-1}).
\end{align*}
\end{small}
Note also that, the event $N$ (the complement of $N^c$) can be represented as
\begin{align*}
N=\cup_{j=1,j\neq i}^L \cup_{k_j=m_j+1}^{n_j} \{T_{k_j}^{(j)}>s\}\cup \cup_{k_i=m_i+2}^{n_i} \{T_{k_i}^{(i)}>s\},
\end{align*}
Thus, we get
\begin{small}
\begin{align*}
&P(M \cap N)\\
&=\sum_{l=1}^{n-\sum_{i=1}^L m_i -1}(-1)^{l+1}\underset{r_1+\ldots +r_L=l}{\sum_{r_1=0}^{n_1-m_1}\ldots \sum_{r_i=0}^{n_i-m_i-1}\ldots \sum_{r_L=0}^{n_L-m_L}}\binom{n_1-m_1}{r_1}\ldots \binom{n_i-m_i-1}{r_i}\ldots\binom{n_L-m_L}{r_L}\\
&\hspace{2cm}\times  \hat{C}(\underbrace{\bar{F}_j(\tau)}_{l_j},\underbrace{\bar{F}_j(s)}_{m_j-l_j+r_j},\underbrace{1}_{n_j-m_j-r_j}, \text{for $j=1,...,L, j\neq i$},  \underbrace{\bar{F}_i(\tau)}_{l_i}, \underbrace{\bar{F}_i(s)}_{m_i-l_i+r_i+1},
\underbrace{1}_{n_i-m_i-r_i-1})\\
\end{align*}
\begin{align*}
&P(M \cap N)-P(M)\\
&=\sum_{l=0}^{n-\sum_{i=1}^L m_i -1}(-1)^{l+1}\underset{r_1+\ldots +r_L=l}{\sum_{r_1=0}^{n_1-m_1}\ldots \sum_{r_i=0}^{n_i-m_i-1}\ldots \sum_{r_L=0}^{n_L-m_L}}\binom{n_1-m_1}{r_1}\ldots \binom{n_i-m_i-1}{r_i}\ldots\binom{n_L-m_L}{r_L}\\
&\times  \hat{C}(\underbrace{\bar{F}_j(\tau)}_{l_j},\underbrace{\bar{F}_j(s)}_{m_j-l_j+r_j},\underbrace{1}_{n_j-m_j-r_j}, \text{for $j=1,...,L, j\neq i$},  \underbrace{\bar{F}_i(\tau)}_{l_i}, \underbrace{\bar{F}_i(s)}_{m_i-l_i+r_i+1},\underbrace{1}_{n_i-m_i-r_i-1})\\
&={\sum_{r_1=0}^{n_1-m_1}\ldots \sum_{r_i=0}^{n_i-m_i-1}\ldots \sum_{r_L=0}^{n_L-m_L}}(-1)^{r_1+\ldots +r_L+1}\binom{n_1-m_1}{r_1}\ldots \binom{n_i-m_i-1}{r_i}\ldots\binom{n_L-m_L}{r_L}\\
&\hspace{2cm}\times  \hat{C}(\underbrace{\bar{F}_j(\tau)}_{l_j},\underbrace{\bar{F}_j(s)}_{m_j-l_j+r_j},\underbrace{1}_{n_j-m_j-r_j}, \text{for $j=1,...,L, j\neq i$},  \underbrace{\bar{F}_i(\tau)}_{l_i}, \underbrace{\bar{F}_i(s)}_{m_i-l_i+r_i+1},\underbrace{1}_{n_i-m_i-r_i-1})
\end{align*}
\end{small}
and
\begin{small}
\begin{align*}
&P(M \cap N \cap K)\\
&=\sum_{l=1}^{n-\sum_{i=1}^L m_i -1}(-1)^{l+1}\underset{r_1+\ldots +r_L=l}{\sum_{r_1=0}^{n_1-m_1}\ldots \sum_{r_i=0}^{n_i-m_i-1}\ldots \sum_{r_L=0}^{n_L-m_L}}
\binom{n_1-m_1}{r_1}\ldots \binom{n_i-m_i-1}{r_i}\ldots\binom{n_L-m_L}{r_L}\\
&\hspace{2cm}\times  \hat{C}(\underbrace{\bar{F}_j(\tau)}_{l_j},\underbrace{\bar{F}_j(s)}_{m_j-l_j+r_j},\underbrace{1}_{n_j-m_j-r_j},\text{for $j=1,...,L, j\neq i$},\underbrace{\bar{F}_i(\tau)}_{l_i}, \underbrace{\bar{F}_i(s)}_{m_i-l_i+r_i},\bar{F}_i{(s+\delta)}, \underbrace{1}_{n_i-m_i-r_i-1})\\
&=\sum_{l=0}^{n-\sum_{i=1}^L m_i -1}(-1)^{l+1}\underset{r_1+\ldots +r_L=l}{\sum_{r_1=0}^{n_1-m_1}\ldots \sum_{r_i=0}^{n_i-m_i-1}\ldots \sum_{r_L=0}^{n_L-m_L}}
\binom{n_1-m_1}{r_1}\ldots \binom{n_i-m_i-1}{r_i}\ldots\binom{n_L-m_L}{r_L}\\
&\times  \hat{C}(\underbrace{\bar{F}_j(\tau)}_{l_j},\underbrace{\bar{F}_j(s)}_{m_j-l_j+r_j},\underbrace{1}_{n_j-m_j-r_j},\text{for $j=1,...,L, j\neq i$},\underbrace{\bar{F}_i(\tau)}_{l_i}, \underbrace{\bar{F}_i(s)}_{m_i-l_i+r_i},\bar{F}_i{(s+\delta)}, \underbrace{1}_{n_i-m_i-r_i-1})+P(M \cap K).
\end{align*}
\end{small}

Similarly,  we have $L=\left\{\cup_{j=1}^L \cup_{k_j=l_j+1}^{m_j} \{T_{k_j}^{(j)}>\tau\}  \right\}.$
 Therefore, we get
\begin{small}
\begin{align*}
P(M \cap L)
&=\sum_{y=1}^{\sum_{j=1}^L (m_j -l_j)}(-1)^{y+1}\underset{d_1+\ldots +d_L=y}{\sum_{d_1=0}^{m_1-l_1}\ldots  \sum_{d_L=0}^{m_L-l_L}}
\binom{m_1-l_1}{d_1}\ldots \binom{m_L-l_L}{d_L}\\
&\times  \hat{C}(\underbrace{\bar{F}_j(\tau)}_{l_j+d_j},\underbrace{\bar{F}_j(s)}_{m_j-l_j-d_j},\underbrace{1}_{n_j-m_j}\text{for $j=1,...,L, j\neq i$}, \underbrace{\bar{F}_i(\tau)}_{l_i+d_i}, \underbrace{\bar{F}_i(s)}_{m_i-l_i-d_i+1}, \underbrace{1}_{n_i-m_i-1}),
\end{align*}
\end{small}
and
\begin{small}
\begin{align*}
&P(M \cap L \cap K)
\\
&=\sum_{l=1}^{\sum_{j=1}^L (m_j-l_j)}(-1)^{l+1}\underset{r_1+\ldots +r_L=l}{\sum_{r_1=0}^{n_1-m_1}\ldots \sum_{r_i=0}^{n_i-m_i-1}\ldots \sum_{r_L=0}^{n_L-m_L}}
\binom{n_1-m_1}{r_1}\ldots \binom{n_i-m_i-1}{r_i}\ldots\binom{n_L-m_L}{r_L}\\
&\hspace{2cm}\times  \hat{C}(\underbrace{\bar{F}_j(\tau)}_{l_j+r_j},\underbrace{\bar{F}_j(s)}_{m_j-l_j-r_j},\underbrace{1}_{n_j-m_j}, \text{for $j=1,...,L, j\neq i$},\underbrace{\bar{F}_i(\tau)}_{l_i+r_i}, \underbrace{\bar{F}_i(s)}_{m_i-l_i-r_i},\bar{F}_i{(s+\delta)}, \underbrace{1}_{n_i-m_i-1})
\end{align*}
\end{small}

Also, after some manipulations, one can verify that $P(M \cap N \cap L)$ and $P(M \cap N \cap L \cap K)$ can be written, respectively, as
\begin{small}
\begin{align*}
&P(M \cap N \cap L)\\
&=\sum_{l=1}^{n-\sum_{j=1}^L m_j -1}(-1)^{l+1}\underset{r_1+\ldots +r_L=l}{\sum_{r_1=0}^{n_1-m_1}\ldots \sum_{r_i=0}^{n_i-m_i-1}\ldots \sum_{r_L=0}^{n_L-m_L}}
\binom{n_1-m_1}{r_1}\ldots \binom{n_i-m_i-1}{r_i}\ldots\binom{n_L-m_L}{r_L}\\
&\times\sum_{y=1}^{\sum_{j=1}^L (m_j-l_j)}(-1)^{y+1}\underset{d_1+\ldots +d_L=y}{\sum_{d_1=0}^{m_1-l_1}\ldots \sum_{d_L=0}^{m_L-l_L}}
\binom{m_1-l_1}{d_1}\ldots \binom{m_L-l_L}{d_L}\\
&\times  \hat{C}(\underbrace{\bar{F}_j(\tau)}_{l_j+d_j},\underbrace{\bar{F}_j(s)}_{m_j-l_j+r_j-d_j},\underbrace{1}_{n_j-m_j-r_j},\text{for $j=1,...,L, j\neq i$},\underbrace{\bar{F}_i(\tau)}_{l_i+d_i}, \underbrace{\bar{F}_i(s)}_{m_i-l_i+r_i+1-d_i},
\underbrace{1}_{n_i- m_i-r_i-1})\\
&=\sum_{l=0}^{n-\sum_{j=1}^L m_j -1}(-1)^{l+1}\underset{r_1+\ldots +r_L=l}{\sum_{r_1=0}^{n_1-m_1}\ldots \sum_{r_i=0}^{n_i-m_i-1}\ldots \sum_{r_L=0}^{n_L-m_L}}
\binom{n_1-m_1}{r_1}\ldots \binom{n_i-m_i-1}{r_i}\ldots\binom{n_L-m_L}{r_L}\\
&\times\sum_{y=1}^{\sum_{j=1}^L (m_j-l_j)}(-1)^{y+1}\underset{d_1+\ldots +d_L=y}{\sum_{d_1=0}^{m_1-l_1}\ldots \sum_{d_L=0}^{m_L-l_L}}
\binom{m_1-l_1}{d_1}\ldots \binom{m_L-l_L}{d_L}\\
&\times  \hat{C}(\underbrace{\bar{F}_j(\tau)}_{l_j+d_j},\underbrace{\bar{F}_j(s)}_{m_j-l_j+r_j-d_j},\underbrace{1}_{n_j-m_j-r_j},\text{for $j=1,...,L, j\neq i$},\underbrace{\bar{F}_i(\tau)}_{l_i+d_i}, \underbrace{\bar{F}_i(s)}_{m_i-l_i+r_i+1-d_i},
\underbrace{1}_{n_i- m_i-r_i-1})+P(M \cap L)
\end{align*}
\begin{align*}
&=\sum_{l=0}^{n-\sum_{j=1}^L m_j -1}(-1)^{l+1}\underset{r_1+\ldots +r_L=l}{\sum_{r_1=0}^{n_1-m_1}\ldots \sum_{r_i=0}^{n_i-m_i-1}\ldots \sum_{r_L=0}^{n_L-m_L}}
\binom{n_1-m_1}{r_1}\ldots \binom{n_i-m_i-1}{r_i}\ldots\binom{n_L-m_L}{r_L}\\
&\times\sum_{y=0}^{\sum_{j=1}^L (m_j-l_j)}(-1)^{y+1}\underset{d_1+\ldots +d_L=y}{\sum_{d_1=0}^{m_1-l_1}\ldots \sum_{d_L=0}^{m_L-l_L}}
\binom{m_1-l_1}{d_1}\ldots \binom{m_L-l_L}{d_L}\\
&\times  \hat{C}(\underbrace{\bar{F}_j(\tau)}_{l_j+d_j},\underbrace{\bar{F}_j(s)}_{m_j-l_j+r_j-d_j},\underbrace{1}_{n_j-m_j-r_j},\text{for $j=1,...,L, j\neq i$},\underbrace{\bar{F}_i(\tau)}_{l_i+d_i}, \underbrace{\bar{F}_i(s)}_{m_i-l_i+r_i+1-d_i},
\underbrace{1}_{n_i- m_i-r_i-1})\\
&+\sum_{l=0}^{n-\sum_{i=1}^L m_i -1}(-1)^{l+1}\underset{r_1+\ldots +r_L=l}{\sum_{r_1=0}^{n_1-m_1}\ldots \sum_{r_i=0}^{n_i-m_i-1}\ldots \sum_{r_L=0}^{n_L-m_L}}
\binom{n_1-m_1}{r_1}\ldots \binom{n_i-m_i-1}{r_i}\ldots\binom{n_L-m_L}{r_L}\\
&\times  \hat{C}(\underbrace{\bar{F}_j(\tau)}_{l_j},\underbrace{\bar{F}_j(s)}_{m_j-l_j+r_j},\underbrace{1}_{n_j-m_j-r_j}, \text{for $j=1,...,L, j\neq i$},  \underbrace{\bar{F}_i(\tau)}_{l_i}, \underbrace{\bar{F}_i(s)}_{m_i-l_i+r_i+1},
\underbrace{1}_{n_i-m_i-r_i})+P(M \cap L).
\end{align*}
\end{small}
and
\begin{small}
\begin{align*}
&P(M \cap N \cap L \cap K)\\
&=\sum_{l=1}^{n-\sum_{j=1}^L m_j -1(t)}(-1)^{l+1}\underset{r_1+\ldots +r_L=l}{\sum_{r_1=0}^{n_1-m_1}\ldots \sum_{r_i=0}^{n_i-m_i-1}\ldots \sum_{r_L=0}^{n_L-m_L}}
\binom{n_1-m_1}{r_1}\ldots \binom{n_i-m_i-1}{r_i}\ldots\binom{n_L-m_L}{r_L}\\
&\times\sum_{y=1}^{\sum_{i=1}^L (m_i-l_i)}(-1)^{y+1}\underset{d_1+\ldots +d_L=y}{\sum_{d_1=0}^{m_1-l_1}\ldots \sum_{d_L=0}^{m_L-l_L}}
\binom{m_1-l_1}{d_1}\ldots \binom{m_L-l_L}{d_L}\\
&\times  \hat{C}(\underbrace{\bar{F}_j(\tau)}_{l_j+d_j},\underbrace{\bar{F}_j(s)}_{m_j-l_j+r_j-d_j},\underbrace{1}_{n_j-m_j-r_j}, \text{for $j=1,...,L, j\neq i$}, \underbrace{\bar{F}_i(\tau)}_{l_i+d_i}, \underbrace{\bar{F}_i(s)}_{m_i-l_i+r_i-d_i},
\bar{F}_i{(s+\delta)},\underbrace{1}_{n_i-m_i-r_i-1})\\
&=\sum_{l=0}^{n-\sum_{j=1}^L m_j -1}(-1)^{l+1}\underset{r_1+\ldots +r_L=l}{\sum_{r_1=0}^{n_1-m_1}\ldots \sum_{r_i=0}^{n_i-m_i-1}\ldots \sum_{r_L=0}^{n_L-m_L}}
\binom{n_1-m_1}{r_1}\ldots \binom{n_i-m_i-1}{r_i}\ldots\binom{n_L-m_L}{r_L}\\
&\times\sum_{y=1}^{\sum_{i=1}^L (m_i-l_i)}(-1)^{y+1}\underset{d_1+\ldots +d_L=y}{\sum_{d_1=0}^{m_1-l_1}\ldots \sum_{d_L=0}^{m_L-l_L}}
\binom{m_1-l_1}{d_1}\ldots \binom{m_L-l_L}{d_L}\\
&\times  \hat{C}(\underbrace{\bar{F}_j}_{l_j+d_j},\underbrace{\bar{F}_j(s)}_{m_j-l_j+r_j-d_j},\underbrace{1}_{n_j-m_j-r_j}, \text{for $j=1,...,L, j\neq i$}, \underbrace{\bar{F}_i(\tau)}_{l_i+d_i}, \underbrace{\bar{F}_i(s)}_{m_i-l_i+r_i-d_i},
\bar{F}_i{(s+\delta)},\underbrace{1}_{n_i-m_i-r_i-1})\\
&+P(M \cap L \cap K)
\end{align*}
\begin{align*}
&=\sum_{l=0}^{n-\sum_{j=1}^L m_j -1}(-1)^{l+1}\underset{r_1+\ldots +r_L=l}{\sum_{r_1=0}^{n_1-m_1}\ldots \sum_{r_i=0}^{n_i-m_i-1}\ldots \sum_{r_L=0}^{n_L-m_L}}
\binom{n_1-m_1}{r_1}\ldots \binom{n_i-m_i-1}{r_i}\ldots\binom{n_L-m_L}{r_L}\\
&\times\sum_{y=0}^{\sum_{i=1}^L (m_i-l_i)}(-1)^{y+1}\underset{d_1+\ldots +d_L=y}{\sum_{d_1=0}^{m_1-l_1}\ldots \sum_{d_L=0}^{m_L-l_L}}
\binom{m_1-l_1}{d_1}\ldots \binom{m_L-l_L}{d_L}\\
&\times  \hat{C}(\underbrace{\bar{F}_j(\tau)}_{l_j+d_j},\underbrace{\bar{F}_j(s)}_{m_j-l_j+r_j-d_j},\underbrace{1}_{n_j-m_j-r_j}, \text{for $j=1,...,L, j\neq i$}, \underbrace{\bar{F}_i(\tau)}_{l_i+d_i}, \underbrace{\bar{F}_i(s)}_{m_i-l_i+r_i-d_i},
\bar{F}_i{(s+\delta)},\underbrace{1}_{n_i-m_i-r_i-1})\\
&+\sum_{l=0}^{n-\sum_{j=1}^L m_j -1}(-1)^{l+1}\underset{r_1+\ldots +r_L=l}{\sum_{r_1=0}^{n_1-m_1}\ldots \sum_{r_i=0}^{n_i-m_i-1}\ldots \sum_{r_L=0}^{n_L-m_L}}
\binom{n_1-m_1}{r_1}\ldots \binom{n_i-m_i-1}{r_i}\ldots\binom{n_L-m_L}{r_L}\\
&\times  \hat{C}(\underbrace{\bar{F}_j(\tau)}_{l_j},\underbrace{\bar{F}_j(s)}_{m_j-l_j+r_j},\underbrace{1}_{n_j-m_j-r_j}, \text{for $j=1,...,L, j\neq i$}, \underbrace{\bar{F}_i(\tau)}_{l_i}, \underbrace{\bar{F}_i(s)}_{m_i-l_i+r_i},
\bar{F}_i{(s+\delta)},\underbrace{1}_{n_i-m_i-r_i-1})\\
&+P(M \cap L \cap K).
\end{align*}
\end{small}

Finally, by replacing the obtained expressions in \eqref{av5} we have

\begin{small}
\begin{align*}
&A_{\mathbf{m}, \mathbf{l}}^{(i)}(s,s+\delta, \tau)=\sum_{l=0}^{n-\sum_{j=1}^L m_j -1}(-1)^{l+1}\underset{r_1+\ldots +r_L=l}{\sum_{r_1=0}^{n_1-m_1}\ldots \sum_{r_i=0}^{n_i-m_i-1}\ldots \sum_{r_L=0}^{n_L-m_L}}
\binom{n_1-m_1}{r_1}\ldots \binom{n_i-m_i-1}{r_i}\\
&\ldots\binom{n_L-m_L}{r_L}\times\sum_{y=0}^{\sum_{j=1}^L (m_j-l_j)}(-1)^{y+1}\underset{d_1+\ldots +d_L=y}{\sum_{d_1=0}^{m_1-l_1}\ldots \sum_{d_L=0}^{m_L-l_L}}
\binom{m_1-l_1}{d_1}\ldots \binom{m_L-l_L}{d_L}\\
&\times  \hat{C}(\underbrace{\bar{F}_j(\tau)}_{l_j+d_j},\underbrace{\bar{F}_j(s)}_{m_j-l_j+r_j-d_j},\underbrace{1}_{n_j-m_j-r_j},\text{for $j=1,...,L, j\neq i$},\underbrace{\bar{F}_i(\tau)}_{l_i+d_i}, \underbrace{\bar{F}_i(s)}_{m_i-l_i+r_i+1-d_i},
\underbrace{1}_{n_i- m_i-r_i-1})\\
&-\sum_{l=0}^{n-\sum_{j=1}^L m_j -1}(-1)^{l+1}\underset{r_1+\ldots +r_L=l}{\sum_{r_1=0}^{n_1-m_1}\ldots \sum_{r_i=0}^{n_i-m_i-1}\ldots \sum_{r_L=0}^{n_L-m_L}}
\binom{n_1-m_1}{r_1}\ldots \binom{n_i-m_i-1}{r_i}\ldots\binom{n_L-m_L}{r_L}\\
&\times\sum_{y=0}^{\sum_{i=1}^L (m_i-l_i)}(-1)^{y+1}\underset{d_1+\ldots +d_L=y}{\sum_{d_1=0}^{m_1-l_1}\ldots \sum_{d_L=0}^{m_L-l_L}}
\binom{m_1-l_1}{d_1}\ldots \binom{m_L-l_L}{d_L}\\
&\times  \hat{C}(\underbrace{\bar{F}_j(\tau)}_{l_j+d_j},\underbrace{\bar{F}_j(s)}_{m_j-l_j+r_j-d_j},\underbrace{1}_{n_j-m_j-r_j}, \text{for $j=1,...,L, j\neq i$}, \underbrace{\bar{F}_i(\tau)}_{l_i+d_i}, \underbrace{\bar{F}_i(s)}_{m_i-l_i+r_i-d_i},
\bar{F}_i{(s+\delta)},\underbrace{1}_{n_i-m_i-r_i-1})
\end{align*}
\begin{align*}
&={\sum_{r_1=0}^{n_1-m_1}\ldots \sum_{r_i=0}^{n_i-m_i-1}\ldots \sum_{r_L=0}^{n_L-m_L}}(-1)^{r_1+\ldots +r_L}
\binom{n_1-m_1}{r_1}\ldots \binom{n_i-m_i-1}{r_i}\ldots\binom{n_L-m_L}{r_L}\\
&\times{\sum_{d_1=0}^{m_1-l_1}\ldots \sum_{d_L=0}^{m_L-l_L}}(-1)^{d_1+\ldots +d_L}
\binom{m_1-l_1}{d_1}\ldots \binom{m_L-l_L}{d_L}\\
&\times \left[  \hat{C}(\underbrace{\bar{F}_j(\tau)}_{l_j+d_j},\underbrace{\bar{F}_j(s)}_{m_j-l_j+r_j-d_j},\underbrace{1}_{n_j-m_j-r_j},\text{for $j=1,...,L, j\neq i$},\underbrace{\bar{F}_i(\tau)}_{l_i+d_i}, \underbrace{\bar{F}_i(s)}_{m_i-l_i+r_i+1-d_i},
\underbrace{1}_{n_i- m_i-r_i-1})\right.\\
&\left.-\hat{C}(\underbrace{\bar{F}_j(\tau)}_{l_j+d_j},\underbrace{\bar{F}_j(s)}_{m_j-l_j+r_j-d_j},\underbrace{1}_{n_j-m_j-r_j}, \text{for $j=1,...,L, j\neq i$}, \underbrace{\bar{F}_i(\tau)}_{l_i+d_i}, \underbrace{\bar{F}_i(s)}_{m_i-l_i+r_i-d_i},
\bar{F}_i{(s+\delta)},\underbrace{1}_{n_i-m_i-r_i-1})\right]\qed
\end{align*}
\end{small}

%

\end{document}